\documentclass[english,11pt]{article}

\usepackage[latin1]{inputenc}
\usepackage{color}

\usepackage{amsfonts}
\usepackage{amsmath}
\usepackage{amssymb}
\usepackage{algorithm2e}
\usepackage{graphicx}

\newcommand{\eps}{\epsilon}
\newcommand{\qed}{$\square$}

\newcommand{\bfv}{{\bf v}}

\newcommand{\dx}{\Delta x}

\newcommand{\NN}{\mathbb{N}}
\newcommand{\RR}{\mathbb{R}}
\newcommand{\CC}{\mathbb{C}}

\newcommand{\ux}{{\underline x}}

\def\eps{\epsilon}

\def\square#1{\vbox{\hrule\hbox{\vrule height#1%
     \kern#1\vrule}\hrule}}          
\def\rectangle#1#2{\vbox{\hrule\hbox{\vrule height#1%
     \kern#2\vrule}\hrule}}
\def\qed{\hfill \vrule height7pt width3pt depth0pt}

\font\tenbb=msbm10
\font\sevenbb=msbm7
\font\fivebb=msbm5

\newfam\bbfam
\scriptscriptfont\bbfam=\fivebb
\textfont\bbfam=\tenbb
\scriptfont\bbfam=\sevenbb

\def\bfv{{\bf v}}

\def\bfN{{\bf N}}

\def\caN{{\cal N}}
\def\caM{{\cal M}}

\def\caF{{\cal F}}

\def\caG{{\cal G}}
\def\caS{{\cal S}}

\def\caM{{\cal M}}

\def\ux{{\underline x}}

\def\uxi{{\underline \xi}}

\def\ut{{\underline t}}

\def\utau{{\underline\tau}}
\def\upi{{\underline \pi}}

\newtheorem{theorem}{Theorem}[section]

\newtheorem{remark}{Remark}[section]

\newtheorem{example}{Example}[section]

\newtheorem{hypothesis}{Hypothesis}[section]

\newtheorem{proposition}[theorem]{Proposition}

\newtheorem{definition}[theorem]{Definition}

\newtheorem{corollary}[theorem]{Corollary}

\newtheorem{lemma}[theorem]{Lemma}

\numberwithin{equation}{section}

\begin{document}

\title{Long Time Boundedness of Planar Jump Discontinuities for 
 Homogeneous Hyperbolic Systems}

\author{Jeffrey Rauch\\ Department of Mathematics, University of Michigan, Ann Arbor}

\date{}

\maketitle

\begin{abstract}   Suppose that 
$L(\partial_t,\partial_x)$ is a homogeneous constant
coefficient strongly hyperbolic partial differential
operator on $\RR^{1+d}_{t,x}$ and  $H$ is a characteristic hyperplane.
Suppose  that in a conic neighborhood of the 
conormal variety of $H$,   the characteristic variety of $L$ is
the graph of a real analytic function $\tau(\xi)$ with
${\rm rank}\,\tau_{\xi\xi}$ identically equal to zero or the maximal possible
value $d-1$.
Suppose that  the source function $f$ is compactly supported
in $t\ge 0$
and piecewise smooth with singularities only
on $H$. Then the solution of $Lu=f$ with $u=0$ for 
$t<0$ is uniformly bounded on $\RR^{1+d}$.
Typically when ${\rm rank}\,\tau_{\xi\xi}\ne 0$ on the conormal variety,
 the sup norm of the  the jump in the gradient of $u$
across $H$ grows linearly with $t$.
\end{abstract}

\section{Introduction}

With 
$A_j\in {\rm Hom} (\CC^k)$ for $1\le j\le d$,
and
$\xi\in \RR^d$, define 
\begin{equation}
\label{eq:L}
L\ :=\ 
\frac{\partial}{\partial t}
\ +\ 
\sum_{j=1}^d
A_j\,
\frac{\partial}{\partial x_j}\,,
\qquad
A(\xi)\,:=\, \sum_j  A_j\, \xi_j\,.
\end{equation}
The operator $L$ maps
$\CC^k$ valued functions to themselves.

\begin{hypothesis} 
\label{hyp:stronghyp}
 The operator $L$ is {\bf strongly hyperbolic}, that is
\begin{equation}
\label{eq:kreiss}
\sup_{\xi\in \RR^d } \
\big\|
e^{i
\sum_{j=1}^d
A_j\,\xi_j
}
\big\|_{{\rm Hom}(\CC^k)}
\ <\
\infty\,.
\end{equation}
\end{hypothesis}

Plancherel's Theorem implies that tempered 
solutions of  $Lu=0$ satisfy for all $t\in \RR$,
$$
\big\| u(t)\big\|_{L^2(\RR^d)}
\ \le\
\sup_{\xi\in \RR^d } \
\big\|
e^{i
\sum_{j=1}^d
A_j\,\xi_j
}
\big\|_{{\rm Hom}(\CC^k)}\
\big\| u(0)\big\|_{L^2(\RR^d)}\,.
$$

The assertions of the abstract are presented and proved 
for the hyperplane $H=\{(t,x)\,:\, x_1=0\}$.  The general case follows by 
a linear change of variable in $\RR^{1+d}$ that preserves
the time $t$.

\begin{hypothesis}
\label{hyp:f}
The source term $f$ is compactly supported in time, 
piecewise smooth, and rapidly decreasing in $x$.  Precisely,
$f\in L^\infty(\RR^{1+d}_{t,x}   )$ 
has support in $\{0\le t\le T<\infty\}$. and the restriction of $f$
to each 
each halfplane
$\{(t,x):\pm\, x_1>0\}$ is infinitely differentiable and for all 
$\alpha, \beta$, $(t,x)^\alpha\partial_{t,x}^\beta f\in L^\infty((t,x):\{(t,x): \pm x_1>0\})$.

\end{hypothesis}

The hyperplane $\{(t,x):\pm\, x_1>0\}$ has conormal variety
$$
\bfN^*\big(\{(t,x):x_1=0\}\big)
\ =\
\{(\tau,\xi)\in \RR^{1+d}\setminus 0\,;\,\tau=0\ {\rm and}\  \xi^\prime=0\}\,.
$$

\begin{definition}
\label{def:smoothvariety}
A point $(\utau, \uxi)\in {\rm Char}(L)$ satisfies the {\bf smooth  variety hypothesis}
when there is a conic neighborhood of $(\utau, \uxi)$ so that in that
neighborhood the characteristic variety is equal to  the graph 
of a real analytic function $\tau(\xi)$.
\end{definition}

The function  $\tau(\xi)$ is homogeneous of degree one.   Therefore
the  second derivative in the radial direction vanishes.   So,
${ \rm rank } \,\tau_{\xi\xi}\le d-1$.

\begin{hypothesis}
\label{hyp:caNsmooth}  The hyperplane $\{(t,x)\,:\,x_1=0\}$  is characteristic
for $L$.  This equivalent to  $\bfN^*(\{(t,x):x_1=0\})
\subset        
{\rm Char}\, L$ and also to  $\det A_1=0$.
The  conormal variety  $\{(\tau,\xi)\ne 0:\tau=0 \ {\rm and}\ \xi^\prime=0\}$
satisfies the smooth variety  hypothesis.
\end{hypothesis}

 .

\begin{definition}
\label{def:ed}
 Denote by $\tau(\xi)$ 
 with $\tau(1,0,\dots,0) =0$  the real analytic
function so that 
the characteristic variety is given by $\{\tau=\tau(\xi)\}$
on a conic neighborhood of $\bfN^*(\{(t,x):x_1=0\})$.
Denote by $\bfv := -\nabla_\xi \tau(1,0,\dots,0)$
the associated group velocity.
\end{definition}

Denote by $u$ the unique solution of the Cauchy problem
\begin{equation}
\label{eq:IVP1}
L\, u \ =\ f\,,
\qquad
u=0\ \ {\rm for}\ \ t\le\ 0\,.
\end{equation}
 It is proved in \cite{courantlax} that $u$ is piecewise smooth
 with singularities only on $\{(t,x) : x_1=0\}$. 
 Moreover, outside the support of $f$ the jumps in $u$ are rigidly transported at the 
 group velocity $\bfv$.  They do not decay as $t\to\infty$
(see Appendix \ref{sec:jumps}). 
 This
 contrasts to the dispersive behavior of solutions that tends
 to spread out  with corresponding amplitudes tending
 to zero.    This paper arose from an attempt, over the last decade
 by Gu\'es and I,
 to quantify the dispersive behavior of internal layers of width
 $\sim \eps$ on time scales $1/\eps$ using ideas from
 \cite{guesrauch1}.  Our computations yielded 
 a leading amplitude  that included a solution $u$ of an equation
 like \eqref{eq:IVP1}.  The consistency of 
 the expansions required that 
 $u\in L^\infty(\RR^{1+d})$.   In the present paper that 
 boundedness is finally proved.
We can now  return to the problem of the diffraction of 
internal layers.

Since the operator $L$ propagates $H^s$ regularity, one might
expect that the solution remains bounded in $H^s(\{\pm x_1>0\})$.
The solution would then
 look like two $H^s$ solutions side by side.
This is not true.
In Appendix \ref{sec:jumps} 
it is shown that when ${\rm rank}\,\tau_{\xi\xi}>0$, 
 the sup norm of 
the jump in $\partial u/\partial x_1$ across $\{x_1=0\}$  
typically grows linearly with time.
 Proofs of piecewise smoothness can also
be found in   \cite{couranthilbert}
Chapter VI Sections 4,5,6, and,  \cite{rauchreed}.
The proof in the appendix benefits from some
simplifications  from analogous advances in geometric optics.

\begin{theorem}
\label{thm:big}
  Assume 
  Hypothesis \ref{hyp:stronghyp}, 
   \ref{hyp:f},
  \ref{hyp:caNsmooth},
  and that with $\tau$ from
  Definition \ref{def:ed} either
  ${\rm rank}\,\tau_{\xi\xi}(1,0,\dots,0)=d-1$ or $\tau_{\xi\xi}$ is identically
  equal to zero.
Then the unique solution $u$ of
\eqref{eq:IVP1}
satisfies $u\in L^\infty(\RR^{1+d})$.
\end{theorem}   

\begin{example}
Consider compactly supported $f$.  Write $f=f_1 +f_2$ where
$f_1$ is smooth and $f_2$ is equal to $f$ cutoff to 
 a small neighborhood $x_1=0$.  The source  $f_2$ is small in $L^2(\RR^{1+d})$.
If the characteristic variety of $L$ has no flat sheets, then
the solution with source $f_1$ tends to zero in $L^\infty(\RR^d)$
as $t\to \infty$ (Corollary 3.3.2 in \cite{rauch2012}).  The solution
with data $f_2$ on the other hand has $L^\infty(\RR^d)$ norm
bounded below and small $L^2(\RR^d)$ norm.    
For large $t$, the singularity
is a tall island of small $L^2$ norm  in a sea of small amplitude
that contains most of the $L^2$ norm.    The transition
from slim peaks to the plateau explains the 
large derivatives for $|t|\gg 1$.
\end{example}

\begin{remark}  \rm 
  Whether the boundedness of $u$
holds without the  rank hypotheses
is an open problem.   One one hand, the flatter is $\lambda$
near $\xi^\prime=0$ the more the hyperbolic problem resembles the
transport operator $\partial_t + \bfv\cdot\partial_x$ and the easier
seem $L^\infty$ estimates.   On the other hand, it is 
conceivable that
there is some focussing phenomenon from the high order Taylor
polynomials of $\tau$
 that leads to amplification
of sup norms.    I conjecture that the transport idea is the 
right one and that $u\in L^\infty(\RR^{1+d}_{t,x})$ holds generally.
\end{remark}

Concerning the pure initial value problem,   Theorem \ref{thm:big} 
implies the following.   Only the last assertion is new.

\begin{corollary}
\label{cor:only}
 Assume 
  Hypothesis \ref{hyp:stronghyp}.
Suppose that over  a neighborhood of $\xi^\prime=0$ the characteristic variety of $L$
consists of $m\le d$ real analytic sheets $\tau=\tau^j(\xi),\, 1\le j\le m$.  Suppose
 that for each $j$,
 ${\rm rank}\, \tau^j_{\xi\xi}$ is identically equal to either $d-1$ or 0 on that 
 neighborhood.
Suppose that $g(x)$ is piecewise smooth with singularities only
on $\{x_1=0\}$ and is rapidly decreasing.   Then the solution of the 
Cauchy problem  $Lw=0$, $w(0)=g$ is piecewise smooth on $\RR^{1+d}$
with singularities only on the $ m$ characteristic hyperplanes passing
through $t=0, x_1=0$.  In addition, $w\in L^\infty(\RR^{1+d})$.
\end{corollary}

\begin{example}  The hypotheses are satisfied for operators with eigenvalues
of constant multiplicity with each sheet of the characteristic variety either
convex, concave, or flat.  {\rm For example, Dirac equations and Maxwell's equations.}
\end{example}

\begin{remark}  \rm
The hypothesis is violated only when over $\xi^\prime=0$ there is either an
eigenvalue crossing or a point of a non flat sheet where ${\rm rank}\, \tau_{\xi\xi}<d-1$.
Each of these occurences is a rare event.  {\rm For example, consider symmetric
hyperbolic $L$.  The space of such operators is parametrized by $d$ hermitian
symmetric matrices, and the hyperplanes by their unit conormals.   This is a subset of 
$\RR^M$ described by real polynomial equations.   Crossings are given by
a vanishing discriminant and low rank $\tau$ by additional real algebraic
equations.  The exceptional set is a finite union of  real
algebraic subsets of positive
codimension in the set of all problems.}
\end{remark}

The proof of Theorem \ref{thm:big} starts by  decomposing $u$ microlocally in $x$ only.

\begin{hypothesis}
\label{hyp:chi}
  $\caN\subset \RR^d_\xi\setminus 0$ is an centro-symmetric, convex, open
conic neighborhood of $\xi^\prime=0$ on which $\tau(\xi)$ is real analytic.
$\chi\in C^\infty(\RR^d\setminus 0)$ satisfies $\chi(\sigma\xi) = \chi(\xi)$ for all
$\sigma\in \RR\setminus 0$, and, the support of $\chi$ meets $S^{d-1}$
on a compact subset of $\caN$.  In addition,  $\chi$ is equal to one on
a conic neighborhood $\caN_1$ of  $\{\xi^\prime=0\}$.
\end{hypothesis}

Fourier multipliers with symbols independent of $x$ are used  throughout.   To a symbol
$\chi(\xi)$ is associeated the operator 
$
\chi(D_x) $,
$D_x := i^{-1} \partial_x$
  defined using the Fourier Traansform $\caF$ by 
$$
\caF\big( \chi(D_x) \, f\big) \ :=\ \chi(\xi)\, \caF f\,.
$$
   Analogous formulas apply to operators
in $D_{t,x}$ using the Fourier transform on $\RR^{1+d}_{t,x}$.

For the  microlocal cutoff $\chi(D_x)$ to a neighborhood of $\xi^\prime =0$, 
write
\begin{equation}
\label{eq:decomposition2}
u \ =\
 \chi(D_x) u\ +\ 
  \big(
  I-\chi(D_x) \big)\,
  u\,.
  \end{equation}
  The second summand is the easy part microlocalized where
  $\xi^\prime\ne 0$  where $f$ is microlocally smooth.
  This term is treated    in Section \ref{sec:xiprimenotzero}.

The term $\chi(D_x)u$ is decomposed using the 
spectral projectors
of $A(\xi)$.  The easy case of the Kreiss Matrix Theorem
(see Appendix 2.I in \cite{rauch2012})  asserts that  
 \eqref{eq:kreiss}
holds
 if and only if the next Conditions {\bf A} and {\bf B} hold.
 \vskip.1cm
 
 Condition {\bf A.}   {\it For all $\xi \in \RR^d$, the eigenvalues of $A(\xi)$
 are real and the 
 eigenspaces  
 span $\CC^k$,}
 $$
 \oplus_{\lambda\in {\rm Spec}\, A(\xi) } \
 E_\lambda(\xi) \ =\ 
 \CC^k,
 \qquad
 E_\lambda(\xi) := \{v\in \CC^k:A(\xi) v = \lambda v\}\ =\ {\rm Ker}\, (A(\xi)-\lambda I )\,.
 $$ 
 For $\lambda\in {\rm Spec}\, A(\xi)$, 
 ${\rm Range}\,
(A(\xi)-\lambda I )\, =\,
 \oplus_{\nu\in {\rm Spec}\, A(\xi)\,\setminus \lambda} \, E_\nu(\xi)$.
 Denote by $\pi_{\lambda}(\xi)$ projection along ${\rm Range}\,
(A(\xi)-\lambda I ) $ onto 
$
{\rm Ker} \,(A(\xi)-\lambda I ) $.
 Then
 \begin{equation}
 \label{eq:spectralrepresentation}
 A(\xi) \ =\ 
 \sum_{\lambda\in {\rm Spec}\, A(\xi)} \lambda\ \pi_\lambda(\xi)
 ,
 \quad
 \sum_{\lambda\in {\rm Spec}\, A(\xi)} \ \pi_\lambda(\xi) \ =\ I,
 \quad\ \ 
 \pi_\lambda(\xi) \  \pi_\nu(\xi)  = 0
 \ \ 
 {\rm for}\ \  
 \lambda\ne \nu.
 \end{equation}

 \vskip.2cm
Condition  {\bf B.} {\it  The function $\xi \mapsto \max_{\lambda\in {\rm Spec}\, A(\xi)}\| \pi_\lambda(\xi)\|$
 is uniformly bounded on $\RR^d$.}

\vskip.2cm

With  $\tau(\xi)$
from Definition \ref{def:ed},
the map $\xi \mapsto\pi_{-\tau(\xi)}$ is real analytic on $\caN$ and
for all $s\in \RR\setminus 0$,
$\tau(s\xi) = \tau(\xi)$.
The eigenvalue identities read
\begin{equation}
\label{eq:ev}
\forall \xi\in \caN,\qquad
\big(\sum_j A_j \xi_j \big)\
\pi_{-\tau(\xi)} \ =\ 
\pi_{-\tau(\xi)} \ \big(\sum_j A_j \xi_j \big)\
=\
-\,\tau(\xi)\
\pi_{-\tau(\xi)}\,.
\end{equation}

The decomposition of $\chi(D_x)\, u$ is
\begin{equation}
\label{eq:decomposition3}
\chi(D_x) \, u
\ =\
 \pi_{-\tau(D_x)}\,
\chi(D_x) \, u
 \ +\
 \big(
 I\ -\
  \pi_{-\tau(D_x)}\big)\chi(D_x) \, u\,.
 \end{equation}
 The second summand is reated 
 in 
 Section \ref{sec:oneminuspi}
 by deriving an equation for it that is 
 microlocally 
 elliptic.
The heavy lifting is the analysis of 
the first summand $\pi_{-\tau(D_x)}  \, \chi(D_x) \, u$.
A scalar hyperbolic equation satisfied by this
 part is derived in
Section
\ref{sec:scalarequation}.
The analysis of that equation microlocally at nonstationary points
is presented in Section \ref{sec:nonstationary}.
In Section  \ref{sec:paraxial} the stationary contributions 
are written as the sum of a paraxial approximation and an error term,
$$
\pi_{-\tau(D_x)}  \, \chi(D_x) \, u
\ =\ 
u_{paraxial} \ +\ 
\big(
u\ -\
u_{paraxial}
\big)\,.
$$
The proof  that 
 the paraxial approximation is bounded uses 
Van de Corput's Lemma.
The proof that the error term is bounded and tends to zero as $t\to \infty$
proceeds by a high/low
frequency decomposition in Section
\ref{sec:stationary}.  
The proof that the low frequency term
is bounded 
 requires an
 inequality of stationary phase
 for test functions with $m$ derivatives with $m>d/2$ but close to $d/2$.
   In Appendix 
\ref{sec:stationaryphase} we present an estimate for 
the limit point case of exactly
$d/2$  derivatives.  The estimate,
possibly new, is
weaker by a factor $|\ln \eps|$ than the standard
estimate.

\section{Analysis away from $\{\xi^\prime =0\}$ }
\label{sec:xiprimenotzero}

\begin{proposition}
\label{prop:away}
  Suppose that $f$ and $u$ are as in 
Theorem \ref{thm:big}.
Suppose that $\beta(\xi) \in C^\infty (\RR^d\setminus 0)$
is homogeneous of degree zero  and vanishes on a conic
neighborhood of $\xi^\prime =0$.
Define
$w\,:=\,\beta(D_x )u$,  
so, $w$ is the unique solution of 
$Lw\, =\, \beta(D_x)  f$ 
that vanishes for $t\le 0$.

\vskip.1cm
{\bf i.}  
For any $s\in \RR$, $
\beta(D_x) \,f \ \in\ 
L^\infty(\RR\,;\, H^s(\RR^d))
$ and is supported in $0\le t\le T$.

\vskip.1cm
{\bf ii.}  $w$ is the unque solution of
$Lw = \beta(D_x)\, f$
that vanishes for $t<0$.
\vskip.1cm
{\bf iii.}  $w\in L^\infty(\RR^{1+d})$.

\end{proposition}

{\bf Proof.}  {\bf i.}  
Since $f$ is piecewise smooth one has
for all $s$
$$
\big\langle
D_{t,x^\prime}
\big\rangle^s f \ \in \
L^2(\RR^{1+d})\,.
$$
For $s\ge 0$ write
$$
\beta(D_x) \,f=
\Big(
 \beta(D_x) \,\big\langle
D_{x^\prime}
\big\rangle^{-s}
\Big)
\Big(\big\langle
D_{x^\prime}
\big\rangle^s
f
\Big)\,.
$$
Since $\langle  \xi\rangle\lesssim\langle \xi^\prime\rangle $
on the support of $\chi$ it follows that
$ \beta(D_x)  \big\langle
D_{x^\prime}
\big\rangle^{-s}$ is bounded from $L^2(\RR^{d})$ to 
$H^s(\RR^d)$.     Therefore $\beta (D_x) \,f \in L^\infty(\RR\,;\,H^s(\RR^{d}))$
and is supported in $0\le t\le T$.

\vskip.1cm

{\bf ii.}   
Follows from {\bf i.}

\vskip.1cm

{\bf iii.}    The Duhamel representation
$$
w(t) \ =\ 
\int_{0}^t
e^{i(t-\sigma) \sum A_jD_j}
\
\beta(D_x) \,f(\sigma)
\
d\sigma
$$
implies that 
$$
\|w(t)\|_{H^s(\RR^d)}
\ \le\
\big\|
e^{\sum iA_j\xi_j}
\big\|_{L^\infty(\RR^d_\xi)}
\
\int_{0}^T \| \beta(D_x) \,f\|_{H^s(\RR^d)}
\ <\ \infty\,.
$$
For $s>d/2$,
this bound uniform for  $t\in \RR$ implies that
$w\in L^\infty(\RR\,;\, L^\infty(\RR^{d}))=L^\infty(\RR^{1+d})$. 
\hfill
\qed
\vskip.2cm

\section{Proof that  $\big(I-\pi_{-\tau(D_x)} )\chi(D_x)u\in L^\infty(\RR^{1+d})$}
\label{sec:oneminuspi}

\begin{proposition}
\label{prop:marcelle}
  With Hypothesis \ref{hyp:chi},
$q:=(I-\pi_{-\tau(D_x)} ) \, \chi(D_x)\,u\in L^\infty(\RR^{1+d})$. 
\end{proposition}

{\bf Proof.}  {\bf I.}  The function  $q$ is supported in $t\ge 0$ and
satisfies two equations, 
$$
Lq = (I-\pi_{-\tau(D_x)} )) \,\chi(D_x)\, f\,,
\qquad
{\rm and},
\qquad
\pi_{-\tau(D_x)} \, q \ =\ 0\,.
$$

From this pair of equations construct a modifed system
as follows.  Choose $\chi_1(\xi)\in C^\infty(\RR^d\setminus 0)$
homogeneous of degree zero, equal to one on a neighborhood
of ${\rm supp}\,\chi \setminus 0$ and 
${\rm supp}\,\chi_1\setminus  0$
is contained in $\caN$ from 
Hypothesis \ref{hyp:chi}. 
Then
\begin{equation}
\label{eq:regularized}
\widetilde L\, q \ =\
 (I-\pi_{-\tau(D_x)} )  \, \chi(D_x)\, f\,,
\qquad
\widetilde L
\ :=\
L \ +\ 
\chi_1(D_x)\ \pi_{-\tau(D_x)} \,
 \partial_{x_1}
\end{equation}
The last summand is everywhere defined since
$\pi_{-\tau(\xi)}$ makes sense  wherever $\chi_1 \ne 0$.

The operator $\widetilde L$ is pseudodifferential in $x$.
In contrast to $L$, the hyperplane $x_1=0$ is noncharacteristic
for $\widetilde L$.   Indeed the symbol
at $\tau=0$, $\xi\in \caN_1$ from Hypothesis 
\ref{hyp:chi}
 is equal to
$$
\widetilde L\big(0, \xi\big) \ =\ i\,A(\xi) + i\, \pi_{-\tau(\xi)}\,\xi_1\,.
$$
This matrix is invertible since $\pi_{-\tau(\xi)}$ is spectral
projection on the kernel of $A(\xi)$. The equation \eqref{eq:regularized}
is therefore elliptic on the wavefront set of the right hand side.

Choose $\delta>0$ so that on $\{|\tau|>\delta\}\times\caN_1$,
$$
\forall (\tau,\xi)\,\in\,
\big\{|\tau|>\delta\big\}\times\caN_1,
\qquad
\big|\det
\big(
\widetilde L\big(0, \xi\big) \ =\ i\,A(\xi) + i\, \pi_{-\tau(\xi)}\,\xi_1
\big)\big|
\ >\ \delta
\,.
$$

\vskip.1cm
{\bf II.}  
For $(\tau,\xi)\in \{|\tau|>\delta\}\times\caN_1$ define
$$
E_1(\tau,\xi) \ :=\ 
\big(
\widetilde L\big(0, \xi\big) \ =\ i\,A(\xi) + i\, \pi_{-\tau(\xi)}\,\xi_1
\big)^{-1}\,.
$$
The symbol $E_1$ is singular at the origin.   Choose 
$\zeta\in C^\infty_0(\RR^{1+d}_{\tau,\xi} ) $, with 
$$
\zeta_1=0 \quad {\rm for}\quad |\tau,\xi|\ <\ 1,
\qquad{\rm and},
\qquad
\zeta_1=1 \quad {\rm for}\quad |\tau,\xi|\ > 2\,.
$$
Choose $\zeta_2\in C^\infty(\RR^{1+d}_{\tau,\xi}\setminus 0)$
homogeneous of degree zero so that the support intersects $S^d$
in a compact subset of $\caN_1$ and so that $\zeta_2=1$
on a conic neighborhood $\caN_2$ of $\{\tau=0\,,\, \xi^\prime=0\}$.
Recall the classical space of symbols of order $m\in \RR$ that are
indepent of $t,x$,
$$
a(\tau,\xi)  \, \in\,
S^m(\RR^{1+d}_{\tau,\xi}
\ \ 
\Longleftrightarrow
\ \ 
\forall \alpha,
\ \exists C<\infty,
\quad
\big|\partial_{\tau,\xi}^\alpha a\big|
\, \le\,
C\,
\langle \tau ,\xi\rangle^{m-|\alpha|},
\quad\ \
\langle \tau,\xi\rangle \ :=\ (1+|\tau,\xi|^2)^{1/2}.
$$
In addition $S^{-\infty}:=\cap_m S^m$.
Define
$$
E_2(\tau,\xi)\ :=\
\zeta_1\, \zeta_2\, E_1(\tau,\xi)\ \ \in \ S^{-1}(\RR^{1+d}_{\tau,\xi } ),
\quad
{\rm so},
\quad
\widetilde L
\ E_2(D_{t,x})\,-\,I\ =\ 
R_1(D_{t,x})\ +\ 
R_2(D_{t,x})\,,
$$
with
\begin{equation}
\label{eq:roger}
R_1\ \in \
S^{-\infty}(\RR^{1+d}_{\tau,\xi})\,,
\quad
{\rm and}\quad
R_2\ \in \ S^0(\RR^{1+d}_{\tau,\xi})\ \ {\rm with}\ \  
R_2(\tau,\xi)\ = \ 0
\quad
{\rm on}
\quad
\{|\tau | <\delta\}\times\caN_1
\,.
\end{equation}

$E_2(D_{t,x})$ is of the  form $K*$ with 
${ singsupp}\,K$ equal to the origin in $\RR^{1+d}_{t,x}$.  Choose $\gamma\in C^\infty_0(\RR^{1+d}_{t,x})$
equal to one on a neighborhood of the origin and supported in 
$| t, x |\le 1$ and define
$$
E(D_{t,x}) \ :=\ (\zeta K)*\,,
\quad
{\rm so},
\quad
E-E_2\ \in \ {\rm Op}(S^{-\infty}(\RR^{1+d}_{\tau,\xi}) ) \,.
$$
With a new $R_1$ satisfying 
\eqref{eq:roger}
one has
$$
\widetilde L
\ E(D_{t,x})\,-\,I\ =\ 
R_1(D_{t,x})\ +\ 
R_2(D_{t,x})\,,
$$
The  Schwartz 
 kernel of $E$ is supported at distance  less than one from the diagonal.  
 Since $\widetilde L$ is local in time, the kernel of the $R_j$ is 
 supported in $\{\big((t,x)\,,\, (s,y)\big)\,:\, |t-s|\le 1\}$.
 
{\bf III.}  Define
$$
q_1 \ :=\ E(D_{t,x}) \,  (I-\pi_{-\tau(D_x)})  \, \chi(D_x)\, f\,.
$$
Since $(I-\pi_{-\tau(D_x)})  \, \chi(D_x)\, f\in L^2(\RR^{1+d})$ and 
$E\in {\rm Op} S^{-1}(\RR^{1+d})$, one has  $q_1\in H^1(\RR^{1+d})$.  The support property of the kernel
of $E$ implies that $q_1$ is 
supported in $-1\le t\le T+1$.  More generally,
for any $\alpha$,
\begin{equation}
\label{eq:inhom}
D_{t,x^\prime}^\alpha \,q_1 
\ =\
E(D_{t,x}  ) \,  (I-\pi_{-\tau(D_x)})  \, \chi(D_x)\, D_{t,x^\prime}^\alpha  f
\ \in \
H^1(\RR^{1+d}   )\,.
\end{equation}
The inclusions \eqref{eq:inhom} 
imply,  by 
a 
Sobolev embedding,
that $q_1\in L^\infty(\RR^{1+d})$.
{\it To complete the proof it suffices to show that 
$q-q_1\in L^\infty(\RR^{1+d})$.}

The function $q_1$ satisfies
$$
\widetilde L q_1 \ =\ (I-\pi_{-\tau(D_x)})  \, \chi(D_x)\, f \ +\ 
g_1 + g_2\,,
\qquad
g_j := 
R_j(D_{t,x})
\
(I-\pi_{-\tau(D_x)})  \, \chi(D_x)\, f \,.
$$
The kernel support property of the $R_j$ imply that the $g_j$
are supported in $-1\le t\le T+1$.
Define $v_j$ to be the solutions of,
\begin{equation}
\label{eq:tildecauchy}
\widetilde L\, v_j \ =\ g_j\,,
\quad
g_j\ =\ 0 \quad{\rm for}\quad t\le -1\,,
\qquad
{\rm so},
\qquad
q\ = \
q_1\ +\ v_1\ +\ v_2\,.
\end{equation}
{\it It suffices  to show that $v_j\in L^\infty(\RR^{1+d}_{t,x})$.}

\vskip.1cm
{\bf IV.}   The Cauchy problem for $\widetilde L$ has existence, uniqueness and estimates entirely
analogous to those of $L$.  Solutions are explicit on the Fourier transform side
and are justified and estimated using the important symbol estimate
\begin{equation}
\label{eq:kreiss2}
\sup_{t\in \RR, \xi\in \RR^d  } 
\big\|
e^{it(
A(\xi) 
\,+\,
\,\chi_1(\xi)\, \pi_{-\tau(\xi)}\, \xi_1)
}
\big\|_{{\rm Hom}(\CC^k)}
=
\sup_{\xi\in \RR^d } 
\big\|
e^{i(
A(\xi) \,+\,
\,\chi_1(\xi)\, \pi_{-\tau(\xi)}\, \xi_1)
}
\big\|_{{\rm Hom}(\CC^k)}
\, <
\infty\,.
\end{equation}
To  prove \eqref{eq:kreiss2}, use \eqref{eq:spectralrepresentation} to write
\begin{align*}
A(\xi)\ +\  \chi_1(\xi)\, \pi_{-\tau(\xi)}\, \xi_1\ =\ 
(-\tau(\xi) + \chi(\xi)\,\xi_1) \,\pi_{-\tau(\xi)}\ +\ 
 \sum_{\lambda\in {\rm Spec}\, A(\xi)\setminus \{-\tau(\xi)\} } \lambda\ \pi_\lambda(\xi)\,.
\end{align*}
With real valued $a_\lambda(\xi)$ this is an expression of the form
\begin{align*}
A(\xi)\ +\  \chi_1(\xi)\, \pi_{-\tau(\xi)}\, \xi_1\ =\ 
 \sum_{\lambda\in {\rm Spec}\, A(\xi)}  a_\lambda(\xi)\ \pi_\lambda(\xi)\,.
\end{align*}
Therefore
$$
e^{i(A(\xi)\ +\  \chi_1(\xi)\, \pi_{-\tau(\xi)}\, \xi_1)}
\ =\
 \sum_{\lambda\in {\rm Spec}\, A(\xi)}  e^{i\,a_\lambda(\xi) }
 \ \pi_\lambda(\xi)\,.
$$
Since $|e^{ia_\lambda}|=1$ this yields
$$
\big\|
e^{i(A(\xi)\ +\  \chi_1(\xi)\, \pi_{-\tau(\xi)}\, \xi_1)}
\big\|_{\CC^k}
\ \le \ 
 \sum_{\lambda\in {\rm Spec}\, A(\xi)}  
 \ 
 \big\|
 \pi_\lambda(\xi)
 \big\|_{\CC^k}
 \,.
$$
Condition {\bf B} from the characterization of strongly hyperbolic operators
implies \eqref{eq:kreiss2}.

\vskip.1cm
{\bf V.  Proof that $v_1\in L^\infty(\RR^{1+d})$.} 
Equation \eqref{eq:kreiss2} 
 imples that for all $s$,  solutions of  $\widetilde Lu=0$ satisfy for all $t\in \RR$,
$$
\big\| u(t)\big\|_{H^s(\RR^d)}
\ \le\
\sup_{\xi\in \RR^d } 
\big\|
e^{i(
A(\xi)
\,+\,
\,\chi_1(\xi)\, \pi_{-\tau(\xi)}\, \xi_1)
}
\big\|_{{\rm Hom}(\CC^k)}\
\big\| u(0)\big\|_{H^s(\RR^d)}\,.
$$

The source term $g_1$ satisfies
$$
g_1\in \cap_s H^s(\RR^{1+d}),
\qquad
{\rm supp}\, g_1 \subset \{-1\le t\le T+1\}\,.
$$
Therefore  for all $s\in \RR$, 
$
g_1
\in L^1_{compact}(\RR\,;\,  H^s(\RR^d))\,,
$ 
Duhamel's formula applied to
\eqref{eq:tildecauchy} 
yields
\begin{align*}
\| v_1(t)\|_{H^s}
\ &\le \
\sup_{\xi\in \RR^d } 
\big\|
e^{i(
A(\xi) \,+\,
\,\chi_1(\xi)\, \pi_{-\tau(\xi)}\, \xi_1)
}
\big\|_{{\rm Hom}(\CC^k)}
\
\int_{-1}^t
\| g_1(\sigma)\|_{H^s(\RR^d)}\ d\sigma
\cr
\ &\le \
\sup_{\xi\in \RR^d } 
\big\|
e^{i(
A(\xi) \,+\,
\,\chi_1(\xi)\, \pi_{-\tau(\xi)}\, \xi_1)
}
\big\|_{{\rm Hom}(\CC^k)}
\
\int_{-1}^{T+1}
\| g_1(\sigma)\|_{H^s(\RR^d)}\ d\sigma
\,.
\end{align*}
Therefore
$
v_1 \ \in \
L^\infty
\big(
\RR\,;\, 
H^s(\RR^d)
\big)\,.
$
Taking  $s>d/2$, yields
$v_1\in L^\infty(\RR^{1+d})$.

\vskip.1cm
{\bf VI.  Proof that $v_2\in L^\infty(\RR^{1+d})$.}   The source term
$g_2\in L^\infty(\RR\,;\, L^2(\RR^d_x))$ is supported in $-1\le t\le T+1$, 
and  is in the image of $R_2(D_{t,x})$.  Therefore,  its space time Fourier
transform vanishes on $\{|\tau|\le \delta\}\times\caN_1$.

On $\RR^{1+d}\setminus (\{|\tau|\le \delta\}\times\caN_1)$,
$\langle \tau,\xi^\prime\rangle^{-2N}$ belongs to 
$S^{-2N}$.  Therefore there is an element $M_{-2N}(\tau,\xi)\in S^{-2N}(\RR^{1+d}_{\tau,\xi} )$
so that $M_{-2N}=\langle \tau,\xi^\prime\rangle^{-N}$ on the support of 
$R_2(\tau,\xi)$.   Therefore on that support,
$$
M_{-2N}(\tau,\xi)\ \Big(1-\tau^2 -\sum_{j\ge 2} \xi_j^2\Big)^N \ =\ 1\,.
$$
The operator version of this identity yields,
\begin{equation}
\label{eq:matthew}
g_2 \ =\ 
R_2(D_{t,x})\
M_{-2N}(D_{t,x})\
(I-\pi_{-\tau(D_x)})  \, \chi(D_x)\
\big(1-\partial_t^2 - \sum_{j\ge 2} \partial_{x_j}^2\big)^N
 f \,.
\end{equation}
In this use that 
$$
\big(1-\partial_t^2 - \sum_{j\ge 2} \partial_{x_j}^2\big)^N
 f \in L^\infty(\RR\,;\, L^2(\RR^d) ) \,,
 \quad
 {\rm with \ support \ in \ \ \ }
 0\le t\le T\,
 $$
 to find that 
 $$
 (I-\pi_{-\tau(D_x)})  \, \chi(D_x)\
\big(1-\partial_t^2 - \sum_{j\ge 2} \partial_{x_j}^2\big)^N
 f
 \ \in \ 
 L^2(\RR^{1+d})\,.
 $$
 Finally $R_2(D_{t,x})\, M_{-2n}(D_{t,x})\in {\rm Op}\, S^{-2N}(\RR^{1+d})$
 so \eqref{eq:matthew} imples that
 $g_2 \in H^{2N}(\RR^{1+d})$.  In addition, $g_2$ is supported in 
 $-1\le t\le T+1$ so $g_2\in L^1(\RR\,;\, H^{2N}(\RR^d))$.   The same
 Duhamel argument as in {\bf V} implies that $v_2\in L^\infty(\RR^{1+d})$.
 Therefore $q = q_1 +v_1 +v_2\in L^\infty(\RR^{1+d})$.  This completes
 the proof of Proposition \ref{prop:marcelle}.
\hfill
\qed

\section{Scalar equation for   $\pi_{-\tau(D_x)} \, \chi(D_x) \, \,u$}
\label{sec:scalarequation}

\subsection{Derivation of the equation}

\begin{lemma}
\label{lem:buzz}
With the notations of Proposition \ref{prop:away}, the
function 
$w:=\pi_{-\tau(D_x)}\, \chi(D_x) \,u $
is characterized as the 
unique solution of the scalar pseudodifferential  initial value 
problem
\begin{equation}
\label{eq:weq}
\Big(
\partial_t
\ -\ 
i\,\tau(D_x)
\Big)
w\ =\ 
\pi_{-\tau(D_x)}  \,
\chi(D_x)
\,f
\,,
\qquad
w=0 \quad {\rm for}
\quad
t<0\,,
\end{equation}
whose spatial Fourier Transform has support in $\caN$.
\end{lemma}

{\bf Proof.}  Define
$v:=\chi(D_x)u$.  Then $v$ is uniquely 
characterized by
\begin{equation}
\label{eq:ger}
L\,v \ =\
\chi(D_x)\,f\,,
\qquad
v = 0 \quad {\rm for}
\quad
t<0\,.
\end{equation}
For any $t$ both the left and right hand side of \eqref{eq:ger} 
have Fourier transforms supported in 
${\rm supp}\,\chi\subset \caN$ so contained in the domain of analyticity
of $\pi_{\tau(\cdot)}$.
Multiply \eqref{eq:ger} by $\pi_{-\tau(D_x)} $  to find
\begin{equation}
\label{eq:isa}
\pi_{-\tau(D_x)} 
\Big(
\partial_t \ +\ 
\sum_j A_j \partial_j\Big)\, \chi(D_x)\, u
\ =\ 
\pi_{-\tau(D_x)}  \,\chi(D_x)\, f\,.
\end{equation}
The symbol of the operator 
$\pi_{-\tau(D_x)} 
\big(
\partial_t \ +\ 
\sum_j A_j \partial_j\big)$ is 
$\pi_{-\tau(\xi)} 
\big(
\partial_t \ +\ 
\sum_j A_j i\xi_j\big)$.
Equaton \eqref{eq:ev} implies that
\begin{align*}
\pi_{-\tau(\xi) }  
\Big(
\partial_t \ +\ 
\sum_j A_j\, i\xi_j\Big)
\ =\ 
\Big(
\partial_t \ +\ 
\sum_j A_j\, i\xi_j\Big)
\pi_{-\tau(\xi)}
\ =\
\Big(
\partial_t
\ -\ 
i\,\tau (\xi)
\Big)
\pi_{-\tau(\xi)} 
\,.
\end{align*}
Therefore
\begin{equation*}
\pi_{-\tau(D_x)} 
\Big(
\partial_t \ +\ 
\sum_j A_j \partial_j\Big)
\ =\ 
\Big(
\partial_t
\ -\ 
i\,\tau(D_x)
\Big)\ 
\pi_{-\tau(D_x)} \,.
\end{equation*}
Injecting this is \eqref{eq:isa}
yields \eqref{eq:weq}
completing the proof.
\hfill
\qed
\vskip.2cm

\subsection{Two simplifications of \eqref{eq:weq} }

\begin{definition}  Denote by $\lambda(\xi)$ the eigenvalue 
$-\tau(\xi)$ of the matrix $A(\xi)$ real analytic on
$\caN$ from Hypothesis
\ref{hyp:chi}.
 The scalar operator in 
 \eqref{eq:weq} is then $\partial_t +i\lambda(D_x)$.
\end{definition}

\subsubsection{Transmission condition and Duhamel}
\label{sec:r}
The symbol
 $\pi_{\lambda(\xi)}\, \chi(\xi)$ is homogeneous of degree
zero and even in $\xi$.  
This implies  that is satisfies the tranmission condition of 
Boutet de Monvel (\cite{boutet}, \cite{H3})
guaranteeing that  $\pi_{\lambda(D_x)}\, \chi(D_x)$  maps piecewise
smooth functions to themselves.   The condition requires that
the Fourier Transform of the distribution
$\pi_{ \lambda (\xi_1, 0,\dots,0)  } \, \chi(\xi_1, 0,\dots,0)$ on $\RR_{\xi_1}$
has smooth extention to each closed half line
$\pm x_1\ge 0$.   The distribution  is identically equal to 
$\pi_{ \lambda (1, 0,\dots,0)  }$ so its Fourier transform is equal to 
a constant times
$\pi_{ \lambda (\xi_1, 0,\dots,0)  } \, \delta(x_1)$.   On each half line
this extends to the smooth function equal to zero.

Choose $\zeta\in C^\infty(\RR^d)$ with $\zeta(\xi)=0$ on $|\xi|<1/2$
and 
$\zeta(\xi)=1$ on $|\xi|>1$.   Write 
$$
\pi_{-\tau(D_x)}  \,
\chi(D_x)
\,f
\ =\
\zeta(D_x)\, 
\pi_{-\tau(D_x)}  \,
\chi(D_x)
\,f
\ +\
(I-\zeta(D_x))\, 
\pi_{-\tau(D_x)}  \,
\chi(D_x)
\,f
\ :=\
F_1 + F_2\,.
$$
The $F_j$ are supported in $0\le t\le T$ 
and have spatial Fourier transform supported in 
${\rm supp}\,\chi$.
Denote by $w_j$ the solutions vanishing for $t\le 0$ with source terms $F_j$
The source $F_2$ has compactly supported Fourier
Transform and
$F_2\in \cap_s\, L^\infty(\RR\,;\, H^s(\RR^d))$.  
Therefore, $w_2\in \cap_s L^\infty(\RR\,;\, H^s(\RR^d) ) \subset
L^\infty(\RR\,;\,L^\infty(\RR^d)) =L^\infty(\RR^{1+d})$.
{\it It  remains to show that $w_1\in L^\infty(\RR^{1+d})$.}

The equation for $w_1$ is
\begin{equation*}
\big(
\partial_t +i\lambda(D_x)
\big) w_1 \ =\ 
F_1\,,
\qquad
w_1=0 \quad {\rm for}\quad
t<0
\end{equation*}
with 
$F_1$ piecwise smooth,  rapidly decreasing,
supported in $0\le t\le T$.
 Duhamel's formula 
 implies that 
 {\it  
 to show that $w_1\in L^\infty(\RR^{1+d})$ it suffices to prove 
  Proposition. \ref{prop:first}.}
 
 \begin{proposition}
 \label{prop:first}
 For any 
$g(x)$ that is piecewise smooth with  singularities in $\{x_1=0\}$,
and rapidly decreasing,
the solution of the initial value problem with
${\rm supp}\, \widehat v(t)\subset    {\rm supp}\,\chi  $,
\begin{equation}
\label{eq:IVP2}
\big(
\partial_t +i\lambda(D_x)
\big) v \ =\ 0\,,
\qquad
v(0)\ =\ 
\chi(D_x)\,g
\end{equation}
satisfies $v\in L^\infty(\RR^{1+d})$.
 \end{proposition}

 \subsubsection{Simpler right hand side with the same jump}

\begin{lemma}  
\label{lem:second}
Suppose that  $\phi(\xi_1)$ is
smooth, even,  identically equal to one for 
$|\xi_1|>2$, and identically  equal to zero for $|\xi_1|\le 1$.
Then Proposition \ref{prop:first} holds provided that for all 
$a\in \caS(\RR^{d-1}_{\xi^\prime} ) $,
\begin{equation}
\label{eq:specialsolution2}
P.V.
\int
\frac
{\chi(\xi)\ 
\phi(\xi_1)\
a(\xi^\prime)}
{\xi_1   }\
e^{i(x\xi + \lambda(\xi)t ) }\
d\xi
\ \in\
L^\infty(\RR^{1+d} _{t,x} )\,.
\end{equation}
\end{lemma}

{\bf Proof of Lemma \ref{lem:second}.}
Denote by $a(x^\prime)\in \caS (\RR^{d-1})$ the jump of $g$ from left to right
at $x^\prime$.
Then the function
$$
\widetilde g\ :=\
e^{-x_1}\
a(x^\prime)
\
{\bf 1}_{x_1>0}
$$
has the same jump as $g$.   It follows that 
for all  $j\le 1$ and $\alpha$
$$
\partial_{x^\prime}^\alpha \, \partial_{x_1}^j \big(g-\widetilde g\big)
\ \in\
L^2(\RR^d)\,.
$$
The solution $v$
of \eqref{eq:IVP2}  
 is the sum of the solution with initial data
$\chi(D_x)\,\widetilde g$ and that with data $\chi(D_x)\big(g-\widetilde g\big)$.   
By inhomogeneous
Sobolev the latter solution belongs to $L^\infty(\RR^{1+d})$.
The former is equal to 
\begin{equation}
\label{eq:specialsolution}
P.V.
\int
\frac
{\chi(\xi)\ 
a(\xi^\prime)}
{1+i\xi_1   }\
e^{i(x\xi + \lambda(\xi)t ) }\
d\xi
\,.
\end{equation}

The solution given by \eqref{eq:specialsolution}
differs from 
\begin{equation}
\label{eq:specialsolution3}
P.V.
\int
\frac
{\chi(\xi)\ 
\phi(\xi_1)\
a(\xi^\prime)}
{1+i\xi_1   }\
e^{i(x\xi + \lambda(\xi)t ) }\
d\xi
\end{equation}
by the inverse Fourier transform of functions uniformly 
in $L^1(\RR^{d})$.  In particular by an element of 
$L^\infty(\RR^{1+d})$.
  It suffices to show that the solution 
given by formula \eqref{eq:specialsolution3}
belongs to 
$L^\infty(\RR^{1+d  } )$.    

Similarly,
\begin{equation*}
P.V.
\int
\frac
{\chi(\xi)\ 
\phi(\xi_1)\
a(\xi^\prime)}
{1+i\xi_1   }\
e^{i(x\xi + \lambda(\xi)t ) }\
d\xi
\ -\
\frac1{i}\
P.V.
\int
\frac
{\chi(\xi)\ 
\phi(\xi_1)\
a(\xi^\prime)}
{\xi_1   }\
e^{i(x\xi + \lambda(\xi)t ) }\
d\xi
\end{equation*}
has spatial Fourier transform belonging to $L^1(\RR^d)$ uniformly in 
$t$.   This completes the proof of Lemma \ref{lem:second}.
\hfill
\qed
\vskip.2cm
The next three sections are devoted to proving
\eqref{eq:specialsolution2}.

 \section{ Nonstationary phase bounds}
 \label{sec:nonstationary}

  \begin{definition}
  \label{def:stationary}
On $\RR^{1+d}_{t,x}\times  \caN$ define the real valued phase
\begin{equation}
\label{eq:defpsi}
\psi(t,x,\xi)\ :=\ x\xi + t\lambda(\xi)\,.
\end{equation}
A point $\ut,\ux,\uxi$ with $0\ne \uxi$ and $\uxi^\prime=0$ is {\bf stationary}
if
\begin{equation}
\label{eq:stationarity}
0\ =\
\nabla_\xi \psi(\ut,\ux,\uxi) \
=\ 
\ux +\ut\lambda_\xi(\uxi) \,.
\end{equation}
\end{definition}

The stationary points with $\uxi^\prime=0$ are exactly those so that
$x= \bfv t$.    
First treat the easier case of nonstationary points  with $\uxi^\prime=0$.  

\begin{hypothesis}
From here on $\phi$ and $a$ are as in Lemma \ref{lem:second}
\end{hypothesis}

\begin{proposition}
\label{prop:nonstationary}
      Suppose that $\ut,\ux,\uxi\in (\RR^{1+d}\setminus 0) \times 
   (\RR^{d}\setminus 0) $ with $\uxi^\prime=0$
   and that 
   $\nabla_\xi\psi(\ut,\ux,\uxi)\ne 0$.   Then there 
   are open conic neighborhoods $\caG$ of $\uxi$ in $\RR^d\setminus 0$ 
   and $\caM$ of 
   $\ut,\ux$ in $\RR^{1+d}\setminus 0$
  so  that if 
 $\gamma(\xi)$ smooth and homogenenous of degree one 
 and whose support intersects $S^{d-1}$ in a  compact 
 subset of $\caG$.  Then for all $\alpha\in \NN^{1+d}$ and $n\in \NN$,
\begin{equation}
\label{eq:arnold}
P.V.
\int
\frac
{\gamma(\xi)\ 
\phi(\xi_1)\
a(\xi^\prime)}
{\xi_1   }\
e^{i(x\xi + \lambda(\xi)t ) }\
d\xi\\ \in\
\langle t,x\rangle^{-n}\
L^\infty\big(
 \caM\big)
\,.
\end{equation}
\end{proposition}

{\bf Proof.}   For $|t|\le 1$ the estimate is easy.
Choose $\caG$ and $\caM$ so that 
   $\lambda$ is smooth on $\caG$ and so that
   there is an $\eta>0$ so that 
      $$
  \big|\nabla \psi\big|\ >\ \eta \,  |t,x| 
  \qquad
  {\rm on}
  \qquad
  \caM     \times
  \caG\,.
   $$

\vskip.1cm
   {\bf I.}  
 Define
on 
$\caM\times\caG$ 
$$
L\ :=\ 
{1
\over
|   \nabla_\xi   \psi   |^2}
\ 
\nabla_\xi \psi
\cdot D_\xi\,,
\qquad
{\rm so},
\qquad
L\, e^{ i\psi}
\ =\ 
e^{ i\psi}
\,.
$$
Using the smoothness of $\lambda$, denote 
by $L^\dagger$  the transposed operator
$$
L^\dagger w \ :=\
D_\xi\bigg(
{
\nabla_\xi\psi
\over
|\nabla_\xi\psi|^2}
\
w\bigg)\,.
$$
The coefficients of $L^\dagger$ are $O( | t,x,\xi|^{-1})$
on the support $\phi(\xi_1)$.  
An integration by parts  yields
$$
P.V.\int
\frac{\gamma(\xi)\
\phi(\xi_1)\
a(\xi^\prime)}
{\xi_1}\
L e^{i\psi(\xi)}
\
d\xi
\ =\ 
P.V.\int
 e^{i\psi(\xi)}\
 L^\dagger
 \Big(
\frac{\gamma(\xi)\
\phi(\xi_1)\
a(\xi^\prime)
}
{\xi_1}
\Big)
\
d\xi
$$
The right hand integral has integrand with $L^1(\RR^d_\xi)$ norm that 
is $O( | t,x|^{-1})$.
Therefore
$$
P.V.\int
\frac{\gamma(\xi)\
\phi(\xi_1)\
a(\xi^\prime)}
{\xi_1}\
e^{i\psi(\xi)}
\
 \frac{1}{\xi_1}\
 d\xi
\ \in\
|t,x|^{-1}\
L^\infty\big(
 \caM\big)
\,.
$$

\vskip.1cm
{\bf II.}   Repeated integration by parts yields for $n\ge 1$,
$$
P.V.\int
\frac{\gamma(\xi)\
\phi(\xi_1)\
a(\xi^\prime)}
{\xi_1}\
L e^{i\psi(\xi)}
\
d\xi
\ =\ 
\int
 e^{i\psi(\xi)}\
 (L^\dagger)^n
 \Big(
\frac{\gamma(\xi)\
\phi(\xi_1)\
a(\xi^\prime)}
{\xi_1}
\Big)
\
d\xi\,.
$$
The operator $(L^\dagger)^n$ in this expression
 has coefficients that are $O(| t,x,\xi|^{-n})$
 on the support of the integrand.
 The right hand integral 
 has integrand with $L^1(\RR^d_\xi)$ norm that 
is $O(| t,x|^{-1})$
implying \eqref{eq:arnold} for $|t|\ge 1$.
\hfill
    \qed

  \begin{remark} \rm
   In addition,  
 for all $\alpha\in \NN^{1+d}$ and $n\in \NN$,
\begin{equation*}
D_{t,x}^\alpha\Big(
P.V.\int
\frac{\gamma(\xi)\
\phi(\xi_1)\
a(\xi^\prime)}
{\xi_1}\
e^{i\psi(\xi)}
\
 d\xi\Big)
\ \in\
\langle t,x\rangle^{-n}\
L^\infty\big(
 \caM\big)
\,.
\end{equation*}
This is not needed below.  The omitted proof follows the strategy  above.
\end{remark}

\section{\bf  Paraxial approximation  for \eqref{eq:specialsolution2}}
\label{sec:paraxial}

Use the homogeneity
$
\lambda(\xi_1,\xi^\prime)
\, =\,
\xi_1
\,\lambda(1,\xi^\prime/\xi_1).
$
Taylor expansion  about $\xi^\prime=0$ yields

$$
\lambda(1,\xi^\prime)\ =\ 
\lambda(1,0) \ +\ 
\xi^\prime\nabla_{\xi^\prime}\lambda(1,0)\ + 
Q(\xi^\prime,\xi^\prime) \ +\ 
{\rm h.o.t}\,,
$$
with
$$
Q(\xi^\prime,\xi^\prime)  \ :=\ 
{1\over 2}
\sum_{2\le i,j\le d-1}
{\partial^2\lambda(1,0)
\over
\partial\xi^\prime_i
\partial\xi^\prime_j
}\
\xi_i\xi_j
\,.
$$
Therefore,
$$
\lambda(1,\xi^\prime/\xi_1)\ =\ 
\lambda(1,0) \ +\ 
\big(\xi^\prime/\xi_1\big)\nabla_{\xi^\prime}\lambda(1,0)\ + 
Q(\xi^\prime/\xi_1,\xi^\prime/\xi_1) \ +\ 
O\big(
|\xi^\prime/\xi_1|^3
\big)\,.
$$
Mulitply by $\xi_1$.  
Using  $\lambda(1,0) = \partial\lambda(1,0)/\partial \xi_1$
in the third line yields,
\begin{equation}
\label{eq:taylor}
\begin{aligned}
\lambda(\xi_1,\xi^\prime)
\ &=\ 
\xi_1\lambda(1,\xi^\prime/\xi_1)
\cr
\ &=\ 
\xi_1
\lambda(1,0) + 
\xi^\prime
\nabla_{\xi^\prime}\lambda(1,0)
+ 
Q(\xi^\prime,\xi^\prime)/\xi_1  + 
O\big(
|\xi^\prime|^3/\xi_1^2
\big)
\cr
\ &=\
\xi\nabla_\xi\lambda(1,0)
\ + \
Q(\xi^\prime,\xi^\prime)/\xi_1 \ +\ 
O\big(
|\xi^\prime|^3/\xi_1^2
\big)
\cr
\ &=\
\bfv\cdot\xi 
\ + \
Q(\xi^\prime,\xi^\prime)/\xi_1 \ +\ 
O\big(
|\xi^\prime|^3/ \xi_1^2
\big)
\,.
\end{aligned}
\end{equation}

 Injecting  in the definition of the solution $u$
 yields the paraxial approximation,
  \begin{equation}
 \label{eq:defparaxial}
 u_{paraxial}
 \ :=\
 P.V. \int
\frac{\chi(\xi)
\
\phi(\xi_1)
\
a(\xi^\prime)
}
{\xi_1}
\
e^{ix\xi}  \
e^{it(\bfv\cdot\xi + Q(\xi^\prime,\xi^\prime)/\xi_1)}
\
d\xi
\,.
\end{equation}

\begin{remark}  The paraxial approximation satisfies 
the differential equation
\begin{equation}
\label{eq:ldpe}
\partial_{\xi_1}\Big(
\partial_t + \bfv\cdot\partial_x\Big)
u_{paraxial}
\ =\ 
Q(D_x^\prime,D_x^\prime)\, u_{paraxial}\,.
\end{equation}
Equation \eqref{eq:ldpe} is classical in  diffractive geometric optics 
with
sources whose spectrum is broad.  For example,
\cite{altermanrauch0}, \cite{altermanrauch2},
\cite{altermanrauch}, \cite{barrailhlannes}, \cite{lannes}.
If $\partial_{\xi_1}$ were replaced by $i$ this would be a Schr\"odinger
equation.   The operators $\partial_{\xi_1}$ and $i$ are both
antiselfadjoint.  The two equations share many properties.
\end{remark}

\begin{theorem}  
\label{thm:paraxialbound}
Suppose that  $0\ne \uxi$,
$\uxi^\prime=0$,
and $t,x,\uxi$ is stationary.  Then the paraxial approximation defined by 
\eqref{eq:defparaxial} satisfies
$
u_{paraxial}\, \in\,
L^\infty(\RR^{1+d}_{t,x})
\,.
$
\end{theorem}

{\bf Proof of Theorem \ref{thm:paraxialbound}.}  The factor $e^{it\bfv\xi}$ induces a translation in $x$ by $\bfv t $.  Thus
it suffices to consider the integral with $\bfv=0$.

The strategy is to integrate $d\xi_1$ with $\xi^\prime$
fixed.  This yields integrals,
\begin{equation}
\label{eq:lim}
\lim_{M\to \infty}\ \int a(\xi^\prime) \Big(\int_{-M}^M
\chi(\xi_1,\xi^\prime ) \
\phi(\xi_1)\  e^{i(x_1\xi_1 + \Lambda/\xi_1)}\
\frac{1}{\xi_1} \ d\xi_1
\Big)\
d\xi^\prime\,,
\qquad
\Lambda\  :=\
t\, Q(\xi^\prime,\xi^\prime)\,.
\end{equation}
Parity in $\xi_1$ shows that the inner integral is equal to
\begin{equation}
\label{eq:M}
\int_0^M
\chi(\xi_1,\xi^\prime ) \
\phi(\xi_1)
\
{
\sin
(x_1\xi_1 + \Lambda/\xi_1)
\over 
\xi_1
}
\
d\xi_1\,.
\end{equation}
Thanks to the rapid decay of $a(\xi^\prime)$, 
{\it to bound the quantity in \eqref{eq:lim}
it suffices to show that there is a constant $C$ so that 
the inegral \eqref{eq:M}
is bounded in absolute value by $C\langle \xi\rangle$
with constant independent of $\Lambda$.}

The next Lemma 
is  the heart of the proof of Theorem \ref{thm:paraxialbound}.
It 
 is proved with a $\Lambda$-dependent high/low 
 frequency decomposition. Van der Corput's Lemma 
treats the low frequency part.

\begin{lemma}
\label{lem:corput}
For all $0<k\in \RR$ there is a constant $C(k)$ so that for 
all real $x,\Lambda$ and all
bounded  intervals $I\subset ]0,\infty[$
one has
$$
\Big|
\int_I
{\sin 
(x\eta+\Lambda / \eta^k)
\over
\eta}
\ 
d\eta\,
\Big|
\ \le \
C(k)\,.
$$
\end{lemma}

{\bf Proof of Lemma \ref{lem:corput}.}   By continuity in $x,\Lambda$ it suffices to treat the
case where $x\ne 0$ and $\Lambda\ne 0$.  Changing the sign
of both $x$ and $\Lambda$ multiplies the integral by $-1$ so it suffices
to consider  $\Lambda>0$.

The domain of integration is divided into two intervals
$I\cap]0,\Lambda^{1/k}[$ and $I\cap]\Lambda^{1/k},\infty[$. 
The first and more interesting is empty when $0<\Lambda<1$.
\vskip.1cm
{\bf I.  Estimate for $I\cap \, ]0,\Lambda^{1/k}[$.}
Change 
variable to
$\xi := e^u$
so,
$$
u=\ln\xi,
\qquad
\xi=e^u,
\qquad
{d\xi
\over
\xi
} \ =\
 du\,.
$$
Then,
$$
x\xi + \Lambda/\xi^k\ =\ \psi(u)
\,,
\qquad
{\rm with}
\qquad
\psi(u)\  :=\ x\,e^u +\Lambda e^{-ku}\,.
$$
The integral is transformed to
$$
\int_J\sin \psi(u)\ du
\,,
\qquad
J\ :=\
\ln I
\ 
\subset \
  \big]-\infty\,,\, \ln(\Lambda^{1/k})   \big[\
  .
$$

With ${}^\prime=d/du$,
$$
\psi^\prime 
= 
x\,e^u \ -\ k\Lambda\,e^{-ku}\,,
\quad
\psi^{\prime\prime} 
=
x\,e^u \ +\ k^2\Lambda\,e^{-ku}\,,
\quad
{\rm and,}
\quad
\psi^{\prime\prime\prime} 
\ =\ 
x\,e^u \ -\ k^3\Lambda\,e^{-ku}
\,.
$$

Since $\Lambda>0$, 
the function $u\mapsto \Lambda e^{-ku}$ is decreasing on $]-\infty,\infty[$.  On $J$ it  
is no smaller than its value at $u=\ln (\Lambda^{1/k})$.  At $u=\ln (\Lambda^{1/k})$,
its value is equal to 1.  

If $x$ and $\Lambda$ have the same  sign,
the summands yielding $\psi^{\prime\prime}$ are both postive.
Thus 
$\psi^{\prime\prime}$ is bounded below by  the 
second summand so on $J$,
  $\psi^{\prime\prime}\ge k^2$.
Van der Corput's Lemma (see \cite{stein}) bounds the integral.  

If $x$ and $\Lambda$ have 
opposite signs the two summands yielding $\psi^{\prime\prime\prime}$ 
are both negative.
The second summand is $\le -k^3$ on $J$.
Therefore   $\psi^{\prime\prime\prime}\le -k^3$ on $J$.
Van der Corput's Lemma bounds the integral.  
Alternatively, $\psi^\prime\le -k$ and $\psi^{\prime\prime}$
is of one sign.  Therefore 
$\psi^\prime$ is monotone in $J$ and again 
Van der Corput's Lemma bounds the integral. 

\vskip.1cm
{\bf II.  Estimate for $I\cap \, ]\Lambda^{1/k},\infty[$.}
First show that 
$$
\int_{I\cap]\Lambda^{1/k},\infty[}
\
{\sin x\xi\over \xi}\ d\xi
\ =\ 
x\,\int_{I\cap]\Lambda^{1/k},\infty[}
\
{\sin x\xi\over x\xi}\ d\xi
$$
is bounded.
Since changing the sign of $x$ multiplies the integral by $-1$
it suffices to treat
$x>0$.  In that case,
 the change of variables
$v=x\xi$ with
$d\xi =x^{-1}dv$ transforms the integral to 
$$
\int_K
{\sin v \over v}\ dv
\,,
\qquad
K\ :=\ x\, \big(
I\cap \, ]\Lambda^{1/k},\infty[
\big)\ \subset\ [0,\infty[\,.
$$
The function $v^{-1}\sin v$ in $v>0$ consists of 
a sequence
of hills of alternating signs and decreasing  areas.
The hill of largest area has area 
equal to $\int_0^\pi (\sin v)/v\,dv$.

The integral over $K$ usually starts and ends with
partial hills both bounded by the largest area.  
The middle is then an alternating decreasing series of hills
whose sum is also bounded by the largest.  Therefore 
the absolute value of the integral over $K$ can be no larger than
$3\int_0^\pi ( \sin v)/v\,dv$.

It then suffices to bound
$$
\Big|
\int_{I\cap]\Lambda^{1/k},\infty[}\ 
{\sin 
(x\xi+\Lambda / \xi^k)
\over
\xi}
\ -\ 
{\sin x\xi\over \xi}
\ 
d\xi
\Big|
$$
Passing the absolute value inside the integral
and using the fact that the derivative of $\sin$
is never larger than one yields the bound
$$
\int_{I\cap]\Lambda^{1/k},\infty[}\ 
{1\over \xi}\
{\Lambda\over \xi^k}\ 
d\xi
\ \le\
\int_{\Lambda^{1/k}}^\infty
\
{\Lambda\over
\xi^{k+1}
}
\ d\xi
\ =\ 
{
\Lambda 
\over
-k}\
\xi^{-k}
\Big|_{\Lambda^{1/k}}^\infty
\ =\ {1\over k}\,.
$$
This completes the proof of  Lemma \ref{lem:corput}.
\hfill
\qed

\vskip.2cm
{\bf End of Proof of 
Theorem \ref{thm:paraxialbound}
. }
It suffices to prove the estimate in italics after 
\eqref{eq:M}.
  Define
$f_\Lambda\in C^1(]0,\infty[ )$ by
$$
f_\Lambda^\prime(\eta) \ =\ 
{\sin 
(x\eta+\Lambda / \eta)
\over
\eta}
\,,
\qquad
f_\Lambda (1)=0\,.
$$
The case $k=1$ of 
Lemma \ref{lem:corput}
implies that $f_\Lambda \in L^\infty(] 0 , \infty [ )$ with bound independent of $\Lambda$.
  Define
  $$
  \qquad
g(\xi) \ =\ g(\xi_1,\xi^\prime) := \chi(\xi ) \
\phi(\xi_1)
\,,
\qquad
h(\xi) \ :=\ \partial g/\partial \xi_1\,.
$$
  For $x_1\ne 0$, the change of variable 
  \begin{equation*}
\label{eq:cov}
\eta=x_1\xi_1,
\qquad
\xi_1= {\eta\over x_1}
\qquad
d\xi_1 \ =\
{d\eta\over x_1}
\end{equation*}
yields
\begin{align*}
\int_0^M
\chi(\xi_1,\xi^\prime ) \
\phi(\xi_1)
\
{
\sin
(x_1\xi_1 + \Lambda/\xi_1)
\over 
\xi_1
}
\
d\xi_1
\ &=\ 
\int_0^{Mx_1}
\chi(\eta/x_1,\xi^\prime)\
\phi(\eta/x_1)
\ 
{\sin 
(x\eta+\Lambda / \eta)
\over
\eta}
\ 
d\eta
\cr
\ &=\
\int_0^{Mx_1}
g\big(\eta/x_1\,,\,\xi^\prime\big)\
f_\Lambda^\prime(\eta)
\ 
d\eta
\,.
\end{align*}

  An integration by parts yields
$$
\ =\ h(M,\xi^\prime)\
f_\Lambda(Mx_1)
\ -\ 
\int_0^{Mx_1}
{d\over d\eta}
\Big(
g(\eta/x_1,\xi^\prime)
\Big)
\
f_\Lambda(\eta)
\ 
d\eta
\,.
$$
The first summand is bounded independent of $M,\Lambda,x_1,\xi^\prime$.
The second summand is equal to 
$$
-
\int_0^{Mx_1}
{1\over x_1}\
\Big(
{d\over d\xi_1}
g(\xi_1,\xi^\prime)
\Big)\Big|_{\xi_1=\eta/x_1}
\ f_\Lambda(\eta)\ d\eta
\ =\ 
-
\int_0^{Mx_1}
{1\over x_1}\
h ( \eta/x_1\,,\xi^\prime)
\ f_\Lambda(\eta)\ d\eta
\,.
$$
The absolute value of this quantity  is bounded by 
$$
\big\|
f_\Lambda\big\|_{L^\infty(\RR_\eta)}\
\Big\|
{1\over x_1}\
h ( \eta/x_1\, , \, \xi^\prime)
\Big\|_{L^1(\RR_\eta)}\,.
$$
Changing back to the variable $\xi_1$ shows that 
$$
\Big\|
{1\over x_1}\
h ( \eta/x_1\, , \, \xi^\prime)
\Big\|_{L^1(\RR_\eta)}
\ =\ 
\int_0^\infty\,
\Big|
h(\xi_1,\xi^\prime)
\Big|
\ d\xi_1\,.
$$
It suffices to show  that 
\begin{equation}
\label{eq:claim}
\int_0^\infty 
\Big|
h(\xi_1,\xi^\prime)
\Big|
\ d\xi_1
\ \lesssim\
\langle \xi^\prime \rangle
\,.
\end{equation}

The function $h$ is continuous and uniformly bounded on $\RR^d$.
Therefore 
it suffices to prove the same bound for the integral with lower limit equal to 2.
   In the range $2\le \xi_1<\infty$, the
function $\phi$ is constant.
Therefore in $|\xi_1|\ge 2$,
  $h$ is homogeneous of degree minus one.
In addition in this range $h(\xi_1,0) =0$.  Therefore
$h=h(\xi_1,\xi^\prime) -h(\xi_1,0)$.   The gradient of $h$ is homogeneous of 
degree minus two.  Estimating the increment in $h$ by a bound for the derivative
times the change in the argument yields
$
|h| \, \lesssim\,
{|\xi^\prime|}/
{|\xi|^2}.
$
This implies  that 
$$
\Big|
\int_2^\infty 
\Big|
h(\xi_1,\xi^\prime)
\Big|
\ d\xi_1
\ \lesssim\
|\xi^\prime|\,.
$$
With the earlier estimate this yields 
\eqref{eq:claim}.
This completes the proof of 
Theorem \ref{thm:paraxialbound}.
\qed
\vskip.2cm

\section{Bound for the stationary contributions}  
\label{sec:stationary}

    It is here that the inequality of stationary
 phase is required for test functions with fractional 
 derivatives.   If one used only integer derivatives it would
 lead to an unnatural lower bound  on the dimension.  The sharp
 limit point inequality of Appendix
 \ref{sec:stationaryphase}
  is more than sufficient.
  
\begin{proposition}
\label{prop:stationary}
Suppose that 
$\uxi^\prime=0$ and that 
  $\lambda$ is stationary at $\ut,\ux,\uxi$. 
    Then there is a conic open neighborhood of $\caG$  
    of  $\uxi$ and
  a conic neighborhood  $\caM$ of $t,x$  so that  
  if $\beta\in C^\infty(\RR^d\setminus 0)$ is homogeneous of 
  degree whose support in $\xi\ne 0$ is contaiined in $\caG$,
$$
P.V.
\int
\frac
{\beta(\xi)\
\phi(\xi_1)\
a(\xi^\prime)}
{\xi_1}
\
e^{ix \xi} 
\
e^{i t\lambda(\xi)}
\
d\xi
\ \in \ 
L^\infty(\caM ) \,.
$$
In addition the left hand side differs from
its paraxial approximation by a term that decays algebraically.  That is,
 there is a $\mu>0$ so that
 $$
PV \hskip-3pt
\int
\frac{\beta(\xi)\,
\phi(\xi_1)\,
a(\xi^\prime)
\,
e^{i(x \xi + t\lambda(\xi) ) }
}{\xi_1}
\,
d\xi
\, -\,
PV \hskip-3pt
\int
\frac{\beta(\xi)\,
\phi(\xi_1)\,
a(\xi^\prime)
\,
e^{i(x \xi + t(\bfv\xi + Q(\xi^\prime,\xi^\prime)/\xi_1))}
}{\xi_1}
\,
d\xi
\,\in\,
t^{-\mu}\, L^\infty(\caM ).
$$
\end{proposition}

\vskip.2cm

{\bf Proof of Propostion \ref{prop:stationary}.}   The result is trivial
when $\tau(\xi)$ is linear in which case the left hand side translates
rigidly at speed $\bfv$ and the paraxial approximation has 
error equal to zero.   Next prove the result when ${\rm rank}\, \lambda_{\xi\xi}$
is equal to $d-1$ on  the support of $\chi$.  The rank hypothesis is only used
in the low frequency bound.

\vskip.1cm
{\bf I.  High frequency bound,
$|\xi_1|\ge t^{1/2 +\delta}$, $0<\delta<1/2$.}

The high frequency bound holds for 
$\caG=\caN$ and  $\caM:=\RR^{1+d}$.
Multiplying \eqref{eq:taylor} by $\xi_1$ yields
$$
\lambda(\xi_1,\xi^\prime)
\ =\ 
\bfv\xi 
\ + \
Q\big(\xi^\prime,\xi^\prime \big)/\xi_1 \ +\ 
O\big( |\xi^\prime|^3/\xi_1^2
\big)
\,.
$$
Therefore
$$
t\lambda(\xi_1,\xi^\prime)
\ -\
t\Big( 
\bfv\xi 
\ + \
Q\big(\xi^\prime,\xi^\prime \big)/\xi_1
\Big) \ =\ 
O\big( \,t|\xi^\prime|^3/\xi_1^2
\big)
\,,
$$
and
$$
e^{  it\lambda(\xi_1,\xi^\prime)  }
\ -\
e^{ i t\big( 
\bfv\xi 
\ + \
Q(\xi^\prime,\xi^\prime )/\xi_1
) }
 \ =\ 
O\big( \,t|\xi^\prime|^3/\xi_1^2
\big)
\,.
$$
Theorem \ref{thm:paraxialbound} shows that 
$$
P.V.\int
\
\beta(\xi)\
e^{i x\xi}
\ 
e^{ i t( 
\bfv\xi 
\ + \
Q(\xi^\prime,\xi^\prime )/\xi_1) }
\
\phi(\xi_1)\
a(\xi^\prime)
\
{
1
\over
\xi_1
}
\
d\xi
\ \in\
L^\infty(\RR^{1+d})
\,.
$$
The difference of this integral and the desired integral is  bounded above
by
$$
C\ 
\int
\
|
\beta(\xi)\
\phi(\xi_1)\
a(\xi^\prime)|
\
{
1
\over
|\xi_1|
}
\
{t|\xi^\prime|^3
\over
|\xi_1^2| }
d\xi
\,.
$$
For $\mu>0$ to be chosen, write the integral as
\begin{equation}
\label{eq:nadal}
t^{-\mu}
\int
\
|
\beta(\xi)\
\phi(\xi_1)\
a(\xi^\prime)|
\
{
1
\over
|\xi_1|
}
\
{t^{1+\mu}|\xi^\prime|^3
\over
\xi_1^2}
d\xi
\,.
\end{equation}
In the high frequency region,
$$
t\ \le \ |\xi_1|^{1/(\delta +1/2)}\,,
\qquad
\frac{t^{1+\mu}}{\xi_1^2}\ \le \ |\xi_1|^{(1+\mu)/(\delta +1/2)-2}\,.
$$ 
Choose
$\mu>0$ so that $(1+\mu)/(\delta +1/2)<2$.
Define 
$0<\nu:=2-(1+\mu)/(\delta +1/2)$.
Then 
 \eqref{eq:nadal} is bounded
by $t^{-\mu}$ times
$$
C\ 
\int
\
|\beta(\xi)
\
\phi(\xi_1)\
a(\xi^\prime)|
\
{
1
\over
|\xi_1|
}
\
{|\xi^\prime|^3
\over
|\xi_1|^\nu }\
d\xi,
$$
Since $a(\xi^\prime)|\xi^\prime|^3\in L^1(\RR^{d-1}_{\xi^\prime})$
and 
$\phi(\xi_1)|\xi_1|^{-1-\nu}\in L^1(\RR_{\xi_1})$,
this integral is absolutely convergent. 
This completes
the high frequency bound.

\vskip.1cm

{\bf II.  Low frequency bound, $|\xi_1|\le t^{1/2 +\delta}$, $0<\delta<1/2$.}

For $t$ fixed, the domain of integration is compact and the integrand is smooth
on that domain.  No principal  value is needed in 
\begin{equation}
\label{eq:scott}
\int_{  |\xi_1| 
  \le t^{1/2+\delta}
  }
  \,
\frac{\beta(\xi)\
\phi(\xi_1)\
a(\xi^\prime)}
{\xi_1}
\
e^{ix \xi} 
\
e^{i t\lambda(\xi)}
\
d\xi\,.
\end{equation}  

Introduce polar coordinates 
\begin{equation}
\label{eq:defomega}
\xi = r\omega,\qquad
|\omega|=1\,.
\end{equation}
Also introduce
\begin{equation}
\label{eq:defw}
w\ :=\ x/t
\qquad
{\rm so},
\qquad
x\,\xi \ =\ t\,r\,w\,\omega\,.
\end{equation}
Stationarity implies that  at $\ut,\ux,\uxi$, $w=\bfv$ is the group velocity at $\uxi$. 
  
  The region of integration
  for the  low frequency region
  lies in 
  $
  0<c<r< t^{1/2+\delta}
  $.
The homogeniety of $x\xi + t\lambda(\xi)$ in $\xi$ implies that 
the integral \eqref{eq:scott}
is equal to 
\begin{equation}
\label{eq:polar}
\int_{c}^{t^{1/2+\delta} }
\Big(
\int_{
S^{d-1}
}
e^{itr(w\omega+\lambda(\omega)) } \
{\beta(r\omega)\
\phi(r\omega_1)\ 
a(r\omega^\prime)
\over
r\omega}
\ 
\,d\omega\Big) \
r^{d-1}
\,dr\,,
\qquad
\xi=r\omega\,.
\end{equation}

 For $w=\bfv$ the phase
 $S^{d-1}\ni \omega\mapsto w\omega+\lambda(\omega)$ is stationary 
 at $\omega=(1,0,\dots,0)$. 
 Choose $\caG$ so small that $\lambda$ is smooth on a neighborhood
 of the closure  of $\caG\cap S^{d-1}$, and  there is an $\eta>0$
 so that on $\caG\cap S^{d-1}$, 
$|\nabla_\xi (w\omega+\lambda(\omega))\ge \eta|\omega^\prime|$.
 
 Since ${\rm rank}\,\lambda_{\xi\xi}(1,0)=d-1$,
 that stationary point is 
nondegenerate on $S^{d-1}$.
  The implicit function theorem implies that for 
$w$ in an open neighborhood $\Omega$ of $\bfv/ \|\bfv\|$ there is a unique  nondegererate
stationary
point in $S^{d-1}$ that lies  close to $(1,0,\dots,0)$.  
Choose $\caM$  so small that for $t,x\in \caM$
this critical point is the only one in $\caG$
and lie in compact subset of $\caG\cap S^{d-1}$.

For
 $w$ in $\caM$ and  $r\in [c_1, c_r t^{1/2+\delta}]$, 
the integral $d\omega$ is a stationary phase integral  with unique
nondegenerate stationary point that is close to $(1,0, \dots , 0)$. 
With the interpolation spaces $Y^m$
defined in 
Appendix 
\ref{sec:stationaryphase},
Theorem \ref{thm:ISP}
yields the 
bound uniform in $w$,
\begin{align*}
\Big|\int_{
S^{d-1}
}
e^{itr(w\omega+\lambda(\omega)) } \ &
{1\over
r\omega_1}\
\phi(r\omega_1)\ 
a(r\omega^\prime)
\ 
d\omega
\Big|
\cr
\ & \lesssim\
|tr|^{-(d-1)/2}\
\ln (1+|tr|)\
\Big\|
{\beta( r\omega)\ \phi(r\omega_1)\ 
a(r\omega^\prime)
\over
r\omega_1}\
\Big\|_{(L^\infty\cap Y^{(d-1)/2} ) (S^{d-1})  }\,.
\end{align*}
The sup norm satisfies, 
\begin{equation}
\label{eq:1overr}
\Big\|
{\beta(r\omega)\ \phi(r\omega_1)\ 
a(r\omega^\prime)
\over
r\omega_1}\
\Big\|_{L^{\infty}(S^{d-1})}
\ \lesssim\
 r^{-1}
 \,.
\end{equation}

$L^1(S^{d-1})$ norms are  smaller.
 First there is the scaling
by $r$ of all $d-1$  variables that
yields,
\begin{equation}
\label{eq:scaling}
\Big\|
{\beta (r\omega )\ \phi(r\omega )\ 
a(r\omega^\prime )
\over
r\omega_1 }\
\Big\|_{L^{1}(S^{d-1})}
\ =\
 r^{-(d-1)}\
\Big\|
{\beta(\zeta)  \ \phi(\zeta) \ 
a(\zeta^\prime)
\over
\zeta_1}\
\Big\|_{L^{1}(rS^{d-1})}
\,.
\end{equation}
Recall that $\zeta=(\zeta_1,\zeta^\prime)$.  Therefore on 
$rS^{d-1}$ with $r$ bounded away from zero, there is a constant 
$c>0$ so that 
$$
c\  {\rm dist}(\zeta, (r,0,0))
\ \le\
|\zeta^\prime |
\ \le\
c^{-1}\
 {\rm dist}(\zeta, (r,0,0))\,.
 $$
The rapid decay of $a$ yields for $\zeta \in rS^{d-1}$,
  $|a(\zeta^\prime)|
  \lesssim 
  {\rm dist}(\zeta, (r,0,0))^{-N}$ for all $N$.
 Therefore the norm on the right of \eqref{eq:scaling} is equal to
$$
 \int_{r\Gamma }
 \Big|
 {\beta (\zeta)\ \phi(\zeta_1)\ 
a(\zeta^\prime)
\over
\zeta_1}\
 \Big|
 \
 d\sigma
 \ \lesssim\ 
  \int_{r\Gamma }
  {
   {\rm dist}(\zeta, \RR(1,0,0))^{-N}
   \over
   r}
   \ d\sigma
 \ \lesssim\ 
 r^{-1} \,.
 $$
 Equation \eqref{eq:scaling}  then yields
$$
\Big\|
{\beta(r\omega)\ \phi( r\omega)\ 
a( r\omega)
\over
r\omega_1 }\
\Big\|_{L^{1}(S^{d-1})}
\ \lesssim\
r^{-d}\,.
$$

   The $r$ in  $\xi = r\omega$ 
appears as a prefactor when one differentiates with  respect to $\omega$
yielding  
\begin{equation}
\begin{aligned}
\label{eq:mminus1}
\Big\|
\partial_\omega^\alpha\,
{\beta(r\omega)\ \phi(r\omega_1)\ 
a(r\omega^\prime)
\over
r\omega_1}\
\Big\|_{L^{1}(S^{d-1})}
\ \lesssim\
r^{-d}\ r^{ |\alpha| }\,.
\end{aligned}
\end{equation}
It follows that  when $d$ is odd so $(d-1)/2$ is an integer, 
\begin{align}
\label{eq:Ybound}
\Big\|
{\beta(r\omega)\ \phi(r\omega_1)\ 
a(r\omega^\prime)
\over
r\omega_1}\
\Big\|_{Y^{(d-1)/2}(S^{d-1})}
\ &\lesssim\
r^{-d}
\
r^{(d-1)/2}
\ =\ 
r^{-(d+1)/2}
\,.
\end{align}
Equation \eqref{eq:Ybound} 
 then holds  by  interpolation when $d$ is even.

The $L^\infty$ contribution 
to the norm in $L^\infty\cap Y^{(d-1)/2}$
is dominant yielding
\begin{align*}
\Big\|
{\beta(r\omega)\ \phi(r\omega_1)\ 
a(r\omega^\prime)
\over
r\omega_1}\
\Big\|_{(L^\infty\cap Y^{(d-1)/2})(S^{d-1}))}
\ &\lesssim\
r^{-1}
\,.
\end{align*}
Therefore, the 
 absolute value of the integral \eqref{eq:polar} 
 is bounded above by
\begin{align*}
\lesssim \
\int_{c}^{t^{1/2+\delta} }
|tr|^{-(d-1)/2}\
&\ln(1+|tr|)
\
r^{-1}
\
r^{d-1}\
dr
\cr
\ &=\ 
t^{-(d-1)/2}
\
\int_{c}^{t^{1/2+\delta} }
\ln(1+|tr|)
\
r^{-1}\
r^{(d-1)/2}\
dr\,.
\end{align*}
Estimate
$$
\ln (1+ |tr|)\ 
\lesssim\
\ln \big(1 + t(t^{1/2+\delta}) \big)
\ \lesssim\
\ln(1+|t|)\,.
$$
The 
 absolute value of the integral \eqref{eq:polar} 
 is therefore
 \begin{align*}
\ \lesssim\
t^{-(d-1)/2}
\
\ln(1+ |t|)
\ 
\int_{c}^{t^{1/2+\delta} }
\
r^{-1}\ r^{(d-1)/2}\
dr
\ \lesssim\
t^{-(d-1)/2}
\
\ln(1+ |t|)
\ 
\big( t^{1/2+\delta}\big)^{(d-1)/2}\,.
\end{align*}
When $\delta<1/2$ this tends to zero as fast as 
a negative  power of $t$  as $t\to \infty$. 
This completes the proof of Proposition \ref{prop:stationary}.
\hfill
\qed

\section{Proofs of  Proposition \ref{prop:first}, Theorem \ref{thm:big}, and Corollary \ref{cor:only}  }

This section combines the preceding results to prove the main Theorems.
  
\vskip.2cm
{\bf Proof of Proposition \ref{prop:first}.}    Thanks to Lemma \ref{lem:second}
it suffices to prove \eqref{eq:specialsolution2}.

   {\bf I. }   {\it If $\uxi = (1,0,\dots,0)$ and $\zeta\in \RR^{1+d}$ is a  unit vector, 
   then there is an open neighborhood $U\subset S^{d}$
of $\zeta$ and $V\subset S^{d-1}$ of  $\uxi$
so that if $\beta\in C^\infty(\RR^d\setminus 0)$ is  homogeneous
of degree zero so that $S^{d-1}\cap{\rm supp}\, \beta\,\subset V$ and 
$\caM\subset \RR^{1+d}$ is the open cone on $U$ then
$$
P.V.\int\,
\frac{\beta(\xi)\
\phi(\xi_1)\
a(\xi^\prime)
}
{\xi_1}
\
e^{ix \xi} 
\
e^{i t\lambda(\xi)}
\
d\xi\ \in\
L^\infty\big(\caM_{\zeta_k}\big)\,.
$$
}
\vskip.1cm

  If $(\zeta, \uxi)$ is nonstationary,
the result follows from Proposition \ref{prop:nonstationary}.
If $(\zeta, \uxi)$ is stationary,
 the result follows from Proposition \ref{prop:stationary}.
This completes the proof of  {\bf I.}

\vskip.1cm

{\bf II.}     For any $\zeta\in S^d$,  choose $U_{\zeta}\subset
S^d$ and $V_\zeta\subset S^{d-1}$ and $M_\zeta$  as in {\bf I.}
Choose a finite subcover $S^{d}\subset \cup_{k=1}^K U_{\zeta_k}$.
Let $W_1:=\cap_{k=1}^K V_{k}$
and define $\caM$ to be the cone on $W_1$, a conic neighborhood of $\xi^\prime=0$.   
Choose an open $W_2\subset S^d$  so that 
$W_2\subset\{\xi^\prime\ne 0\}$ and 
$W_1,W_2$ cover $S^{d-1}$.
Choose
a smooth partiton of unity $1 =\psi_1+\psi_2$ on $S^{d-1}$ with 
${\rm supp}\, \psi_j\subset W_{j}$.   Define $\gamma_j(\xi):=
\psi_j(\xi/|\xi|)$.    The result from {\bf I} implies that 
\begin{equation}
\label{eq:primo}
P.V.\int\,
\frac{
\chi(\xi)\
\gamma_1(\xi) \
\phi(\xi_1)\
a(\xi^\prime)
}
{\xi_1}
\
e^{ix \xi} 
\
e^{i t\lambda(\xi)}
\
d\xi\ \in\
L^\infty\big(  \caM_{\zeta_k}  
 \big)\,,
\ \ 1\le k\le K\,.
\end{equation}
Since the $M_{\zeta_k}$ cover $\RR^{1+d}_{t,x}\setminus 0$
this implies that 
\begin{equation}
\label{eq:primo2}
P.V.\int\,
\frac{
\chi(\xi)\
\gamma_1(\xi) \
\phi(\xi_1)\
a(\xi^\prime)
}
{\xi_1}
\
e^{ix \xi} 
\
e^{i t\lambda(\xi)}
\
d\xi\ \in\
L^\infty\big( \RR^{1+d})
 \big)\,.
\end{equation}

Proposition \ref{prop:nonstationary}
implies that
\begin{equation}
\label{eq:secondo}
P.V.\int\,
\frac{
\chi(\xi)\
\gamma_2(\xi) \
\phi(\xi_1)\
a(\xi^\prime)
}
{\xi_1}
\
e^{ix \xi} 
\
e^{i t\lambda(\xi)}
\
d\xi\ \in\
L^\infty\big(\RR^{1+d}\big)\,.
\end{equation}
Since $\gamma_1 + \gamma_2 =1$,
addiing \eqref{eq:primo} and \eqref{eq:secondo} proves  
\eqref{eq:specialsolution2}.
\hfill 
\qed
\vskip.2cm

{\bf Proof of Theorem \ref{thm:big}.}   Begin with \eqref{eq:decomposition2}.
That $\big(
  I-\chi(D_x) \big)\,
  u\in L^\infty(\RR^{1+d})$ is proved in Proposition 
  \ref{prop:away}.  
  
  The remaining part $\chi(D_x) u$ is decomposed in 
  \eqref{eq:decomposition3}.   That the  summand
  $(I-\pi_{-\tau(D_x)} ) \,\chi(D_x)\, f$ belongs to 
  $L^\infty(\RR^{1+d})$ is proved in 
  Proposition \ref{prop:marcelle}.

  An initial value problem satisfied by the other  summand
  $\pi_{-\tau(D_x)}\chi(D_x)u$ is derived in Lemma \ref{lem:buzz}.  
  In Subsection \ref{sec:r} it is proved that  the boundedness
  of that summand is a consequence of Proposition 
  \ref{prop:first},
  that  has
  just been proved.  
 This completes the proof of Theorem \ref{thm:big}.
\hfill
\qed
\vskip.2cm

{\bf  Proof of Corollary \ref{cor:only}.}
All but the last assertion are proved by constructing such
a solution as a sum of progressing waves as in \cite{couranthilbert}.
To prove boundedness choose a smooth cutoff function $\rho(t)$ vanishing for 
$t<1$ and equal to 1 for $t\ge 2$.   Define $u:=\rho(t)\, u$.
Then $Lu=F$ with  $F$ piecewise 
smooth,   supported in $1\le t\le 2$, and, rapidly decreasing.   The singularites
are on the parts of the $m$ characteristic hyperplanes that lie in $1\le t\le 2$.
Write the source term as $F=\sum_{j=1}^m F_j$ with $F_j$
carrying the singularities on the $j^{\rm th}$.  Then $u=\sum u_j$ where
$u_j$ is the solution to $Lu_j=F_j$ that vanishes in $t\le 0$. 
Apply   Theorem
\ref{thm:big} adapted to sources 
 with singularities on the $j^{\rm th}$
hyperplane  to  complete the proof.
\hfill
\qed

 \appendix

\section{Propagation of jumps across flat discontinuities}
\label{sec:jumps}

Recall that  $\{x_1=0\}$ is characteristic, that is 
$\det A_1=0$.   The spectral projection associated to
$\xi=(1,0,\dots,0)$ and  the 
eigenvalue zero
 has been
denoted
$\pi_0(1,0,\dots,0)$.

\begin{definition}
\label{def:pi} Denote by
$\upi$ the projector 
$\pi_0(1,0,\dots,0)$.
 \end{definition}

\begin{definition}
\label{def:J}
For a piecewise smooth $f(t,x)$ with singularities only on $\{x_1=0\}$
denote by $J_f^n(t,x^\prime)$ the  jump 
$$
J_f^n(t,x^\prime) \ :=\
\frac{\partial^n f}{\partial x_1^n}
\big(t\,,\,
0+\,,\,x^\prime\big)
\ -\
\frac{\partial^n f}{\partial x_1^n}
\big(t\,,\,
0-\,,\,x^\prime\big)\,.
$$
\end{definition}

Then for all $M\ge 0$
$$
f\ -\
\sum_{n=0}^M
J_f^n\
\frac{x_1^n}
{n!}
\
{\bf 1}_{\RR_+}
\ \in\
C^M(\RR^{1+d})\,.
$$
This relation is abbreviated as
$$
f \ =\
\sum_{n=0}^\infty
J_f^n\
\frac{x_1^n}
{n!}\
{\bf 1}_{\RR_+}
\ +\ C^\infty(\RR^{1+d})\,.
$$

Given  such a  source,  Courant and Lax \cite{courantlax}  construct a piecewise smooth
solution $u$ to
\begin{equation}
\label{eq:Luf}
L u \ -\ f\ \in \ C^\infty(\RR^{1+d})\,.
\end{equation}
Direct computation yields
$$
A_1 \,
\frac{\partial u}
{\partial x_1} \ =\ 
A_1\,
 J_u^0 \,
 \delta(x_1) 
\ +\
\sum_{n=0}^\infty
A_1\,
 J_u^{n+1}\, 
\frac{x_1^n}
{n!}
\
{\bf 1}_{\RR_+}
\ +\ C^\infty(\RR^{1+d})\,.
$$

\begin{definition}
\label{def:Ltan}
Denote the tangential part of $L$ by 
$$
L_{tan}\big(\partial_t,\partial_{x^\prime}\big)
\ :=\
\frac{\partial}{\partial t}
\ +\ 
\sum_{j=2}^d
A_j\, 
\frac{\partial}{\partial x_j}\,,
\qquad
{\rm so},
\qquad
L\ =\
L_{tan}+A_1\, \frac{\partial}{\partial x_1}\,.
$$
\end{definition}

Then $L_{tan}$ maps piecewise smooth functions to 
piecewise smooth  functions and 
$$
L_{tan} u \ =\
\sum_{n=0}^\infty
(L_{tan}J_u^n)\
\frac{x_1^n}
{n!}\
{\bf 1}_{\RR_+}
\ +\ C^\infty(\RR^{1+d})\,.
$$

Therefore, in order that \eqref{eq:Luf} be satisfied it is necessary and sufficent that
\begin{equation}
\label{eq:Jeq}
A_1 J_u^0 \ =\ 0,
\end{equation}
and
\begin{equation}
\label{eq:Jeqbis}
\forall n\ge 0,
\qquad
A_1 J_u^{n+1} \ +\ 
L_{tan} \,J_u^n
\ =\ J_f^n \,.
\end{equation}

\begin{definition}
\label{def:Q}
   Define a partial inverse $Q$ to $A_1$ by 
   \begin{equation}
   \label{eq:defQ}
A_1\,\upi \ =\ 0,
\qquad
Q\,A_1\,(I-\upi) \ =\ I-\upi\,.
\end{equation}
\end{definition}

\begin{proposition}
\label{prop:jumps}
Suppose that $f(t,x)$ is 
piecewise smooth with   singularities  only on $\{x_1=0\}$ and 
$f=0$ for $t\le 0$.
Then there are uniquely determined  jumps  $J_u^n\in C^\infty\big( \RR_{t,x^\prime}^d  \big) $
satisfying  the sequence of equations
\eqref{eq:Jeq} and \eqref{eq:Jeqbis}. 
\end{proposition}

{\bf Proof.}  The key to deciphering the equations is to 
observe that  equations \eqref{eq:Jeqbis} hold if and  only if their
projections by $\upi$ and $I-\upi$ hold.   The $\upi$ projection eliminates
the $A_1$ term.  Thus Equations \eqref{eq:Jeqbis}  hold if and only if
\begin{equation}
\label{eq:two}
\forall n\ge 0,
\qquad
\upi \ L_{tan} \ J_u^n =\ \upi \ J_f^n,
\end{equation}
and
\begin{equation}
\label{eq:twobis}
\forall n\ge 0,
\qquad
(I-\upi) A_1 J_u^{n+1} \ =\ 
(I-\upi)
\Big(
J_f^n  \ -\
L_{tan} J_u^n
\Big)
\,.
\end{equation}
Equation
\eqref{eq:two}
 is rewritten by writing $J_u^n
 =
 \upi J_u^n + (I-\upi)J_u^n$
to find
\begin{equation}
\label{eq:first}
\forall n\ge 0,
\qquad
\upi \ L_{tan}\ \upi\, J_u^n \ =\ 
\upi J^n_f \ -\ \upi
 \ L_{tan}\ (I-\upi)J_u^n\,.
\end{equation}
Equation  \eqref{eq:twobis}
is between vectors in the rangle of $I-\upi$.   The equation is 
equivalent to the same equation multiplied by $Q$.   Using
$(1-\upi) A_1 =A_1 (1-\upi)$, \eqref{eq:twobis} is
equivalent to
\begin{equation}
\label{eq:ell}
\forall n\ge 0,
\qquad
(I-\upi)\,J_u^{n+1}  \ =\ Q\Big(
J_f^n  \ -\
L_{tan} J_u^n
\Big)\,.
\end{equation}
Summarizing, {\it  the jumps satisfy
\eqref{eq:Jeq} and \eqref{eq:Jeqbis} if an only 
they satisfy \eqref{eq:Jeq}, \eqref{eq:first}.
and 
\eqref{eq:ell}.
}

In the next disussion, $\eqref{eq:first}_n$ means the case $n$
of equation \eqref{eq:first}.
The jump $J_u^0$ is first determined uniquely
by \eqref{eq:Jeq} and $\eqref{eq:first}_0$.
The  jumps $J_u^n$ with $n\ge 1$
are determined from  $\eqref{eq:first}_n$
and  $\eqref{eq:ell}_{n-1}$.      Conversely, with these
determinations the equations 
$\eqref{eq:Jeq}_{n-1}$,  $\eqref{eq:first}_n$,
and 
\eqref{eq:ell} are  satisfied.
 
First consider $J_u^0$.    Equation \eqref{eq:Jeq} implies
that $(I-\upi )J_u^0=0$.  The case $n=0$ of \eqref{eq:first}
reads $\upi L_{tan}\upi \, J_u^0 \ =\ \upi J_f^0$. 
Since $\upi  A_1=0$ one has $\upi L_{tan} = \upi L$.    The smoothness of 
$\lambda$ on a neighborhood of $(1,0,\dots\,,0)$ implies 
the key transport identity of geometric optics 
(Proposition 5.4.1 in  \cite{rauch2012}),
\begin{equation}
\label{eq:magic}
\upi L_{tan}\upi \ =\ \upi L\upi  \ =\ 
\partial_t \ + \ \bfv\cdot \partial_x\, ,
\qquad
\bfv  \ :=\ \nabla_\xi\lambda(1,0,\dots\,, 0)\,.
\end{equation}
Then $J_u^0=\upi J_u^0$ is determined 
by 
\begin{equation}
\label{eq:eddie}
\Big(
\partial_t \ + \ \bfv\cdot \partial_x
\Big)
 \upi J_u^0 \ =\ 
\upi J^n_f,
\qquad
J_u^0\ =\ 0 \quad 
{\rm for}
\quad
t\le 0\,.
\end{equation}
Conversely,  this equation together with  \eqref{eq:Jeq} imply
 the case $n=0$ of 
\eqref{eq:first} are satisfied.

For $n\ge 1$, the jump $J_u^n$ is uniquely determined
 $\eqref{eq:ell}_{n-1}$ and 
$\eqref{eq:first}_n$.
Begin by replacing  $(I-\upi)J^n_u$  in  of $\eqref{eq:first}_n$ using 
$\eqref{eq:ell}_{n-1}$.  Then using
\eqref{eq:magic} yields
\begin{equation}
\label{eq:first2}
\forall n\ge 1,
\quad
\Big(
\partial_t \, + \, \bfv\cdot \partial_x
\Big)
 \upi J_u^n \, =\, 
\upi J^n_f \ -\ \upi
 \, L_{tan}\,
 Q\,
 L_{tan}\,
J_u^{n-1}\,,
\quad
J_u^n=0 \quad
{\rm for}
\quad
t\le 0\,.
\end{equation}
Solving this detemines $\upi J_u^n$ and 
$\eqref{eq:ell}_{n-1}$ determines $(1-\upi)J_u^n$.
Conversely when $\eqref{eq:first2}_n$ together with 
$\eqref{eq:ell}_{n-1}$ hold, one recovers $\eqref{eq:first}_n$.
\hfill
\qed
\vskip.2cm

\begin{theorem}
Suppose that $f(t,x)$ is a piecewise smooth  function on $\RR^{1+d} $ with
singularities only on $\{x_1=0\}$ and that $f=0$  for $t\le 0$.
Then the unique solution of $Lu=f$ that vanishes for $t\le 0$
is also piecewise smooth with singularities only on $\{x_1=0\}$.
\end{theorem}

{\bf Proof.}  Determine jumps $J_u^n(t,x^\prime)$ using Proposition \ref{prop:jumps}.
Choose a piecewise smooth $v$ vanishing in $t\le 0$ whose jumps are
equal to the functions $J_u^n$.
Then for all $M$, one has
$
L v \, -\, f \, \in\, C^M(\RR^{1+d}).$
Define $g:= Lv-f$.  Then $g\in C^\infty(\RR^{1+d})$ and vanishes for $t\le 0$.
Define $w\in C^\infty(\RR^{1+d})$ to be the solution vanishing for $t\le 0$
to $Lw = -g$.   Then $u:= v+ w$  is piecewise smooth
and 
$$
L u \ =\ Lv \ +\ Lw
\ =\
(g+f) -g \ = \ f\,.
$$
Therefore $u=v+w$ is the unique solution of the Initial value problem.
Since $u$ has
  the desired properties, this  completes the proof.
\hfill
\qed
\vskip.2cm

\begin{corollary}
If $f$ is compactly supported then $J_u^0\in L^\infty(\RR^d_{t,x^\prime} ) $.\end{corollary}

{\bf Proof.}   Choose $R>0$ so that $|t|+|x|>R\Rightarrow f=0$.  The recipe
for $J_u^0$ implies that 
for $|t|>R$
$$
\Big(\partial_t + \bfv\cdot \partial_x \Big)
J_u^0\ =\ 0.
$$
Therefore 
$$
\big\|
J_u^0
\big\|_{L^\infty(\RR^{1+(d-1)} ) }
\ \le\
\big\|
J_u^0
\big\|_{L^\infty(\{|t|\le R\} \times\RR_{x^\prime}^{d-1}) } \,.
 $$
In $|t|\le R$, $J_u^0$ is smooth and compactly supported so bounded.
This shows that the right hand side is finite completing the proof.
\hfill
\qed
\vskip.2cm

\begin{example}   Tangential derivatives $\partial_{t,x^\prime}^\alpha u$
satisfy an equation entirely analogous to that satisfied by $u$.
It follows that the jumps in these derivatives belong to $L^\infty(\RR^d_{t,x^\prime})$.
\end{example}

 Nontangential derivatives need not be bounded.  
 Beyond the support of $f$,
 $$
 \Big(
\partial_t \, + \, \bfv\cdot \partial_x
\Big)
 \upi J_u^1 \, =\, 
 -\ \upi
 \, L_{tan}\,
 Q\,
 L_{tan}\,
J_u^{0}\,.
$$
In this range,
$$
 \Big(
\partial_t \, + \, \bfv\cdot \partial_x
\Big)
\Big(
 \upi
 \, L_{tan}\,
 Q\,
 L_{tan}\,
J_u^{0}
\Big)
\ =\ 0\,.
$$
Therefore {\it either $ \upi
 \, L_{tan}\,
 Q\,
 L_{tan}\,
J_u^{0}
$
is  identically equal to zero beyond the support of $f$,
or $J_u^1$ grows linearly in time.}

Since $\upi J^0_u = J^0_u$ and  $\upi A_1 =A_1 \upi =0$,
$$
 \upi
 \, L_{tan}\,
 Q\,
 L_{tan}\,
J_u^{0}
\ =\ 
 \upi
 \, L_{tan}\,
 Q\,
 L_{tan}\,
 \upi\
J_u^{0}
\ =\
 \upi
 \, L\,
 Q\,
 L\,
 \upi\
J_u^{0}
\,.
$$

The fundamental identity of diffractive geometric optics
\cite{djmr1996} Proposition 3.2,  and, \cite{jmr1998}  reads   
$$
-\
\upi
 \, L\,
 Q\,
 L\,
 \upi
  \ =\ 
\  \frac{1 }2\
\upi\
  \sum
  \frac{\partial^2\lambda(1,0,\dots,0)}
  {\partial \xi_\mu\,
\partial\xi_\nu}\
 \frac{\partial^2}
  {\partial x_\mu\,
\partial x_\nu}\ :=\
P(\partial_{x^\prime})
$$
Therefore  beyond the support of $f$ 
\begin{equation}
\label{eq:JPf}
 \Big(
\partial_t \, + \, \bfv\cdot \partial_x
\Big)
 \upi J_u^1 \, =\, 
 P(\partial_{x^\prime})
 \,
 J^0_f
 \ =\ 
 J^0_{P(\partial_{x^\prime} )f}
 \,.
 \end{equation}
The right hand side 
is constant on integral curves of $\partial_t+\bfv\cdot\partial_x$.
When that constant is not zero, 
$J^1_u$ grows linearly 
along the integral curve.

By hpothesis the matrix of second derivatives of $\lambda$ has rank
$d-1>0$ so the operator $P$ is not identically equal to zero.

Equation \eqref{eq:JPf} asks one to integrate the compactly supported 
 jump in $Pf$
along integral curves of $\partial_t + \bfv\cdot \partial_x$.
The constant in the preceding paragraph vanishes if and only if
 the integration yields answer zero..
 The value zero is a rare occurence.
 {\it For generic $f$, $\upi J^0_u$  is nonzero
and constant on 
almost all (in the sense of Lebesgue measure) integral curves that touch the support of 
$J^0_f$.   For those $f$ and  integral curves, $J^1_u$ grows linearly 
in time.}

 \section{Limit case stationary phase inequality}
  \label{sec:stationaryphase}
  
\subsection{Non stationary phase lemmas}

 \begin{lemma}[Lemma of Nonstationary Phase]
 \label{lem:nonsta1}
 Suppose that $\Omega\subset\RR^d$ is 
open,
 $m\in \NN$, and $C_1>1$.  Then there is a constant $C>0$ so that
for all $f\in C^m_0(\Omega)$, and
$\phi\in C^m(\Omega\,;\,\RR)$ satisfying
$$
\forall |\alpha|\le m+1\,,
\ \ \|\partial^\alpha\phi\|_{L^\infty(\Omega)}\le C_1\,,
\quad
{\rm and},
\quad
\forall x\in \Omega\,,\ \
C_1^{-1}\ \le |\nabla_x\phi|\ \le \ C_1\,,
$$
one has,
\begin{equation}
\label{eq:nonsta}
\Big|
\int \ e^{i\phi (x)/\eps}\ 
f(x)\ dx
\Big|
\ \le \
C\, \eps^m\,
\sum_{|\alpha|\le m}\,
\|\partial^\alpha f\|_{L^1(\RR^d)}\,.
\end{equation}
\end{lemma}

{\bf Proof.}   Introduce the differential operator
of order one with smooth coefficeints
$$
L\ :=\ 
\frac{\nabla \phi}
{  i  |\nabla \phi|}
\ \cdot\
\nabla_x\,,
\qquad
{\rm so},
\qquad
 L e^{i\phi/\eps} \ =\
\
\eps^{-1}\
e^{i\phi/\eps} \,.
$$
Write
$$
\int \ e^{i\phi(x) /\eps}\ 
f(x)\ dx
\ =\
\eps^m\
\int L^m(e^{i\phi/\eps})\
f(x) \dx\
\ =\ 
\eps^m\
\int 
e^{i\phi/\eps}
\ 
\big(
L^\dagger\big)^m f\ dx,
$$
where $L^\dagger$ denotes the transposed operator.
Since  
$$
\big\|
(L^\dagger)^m f
\big\|_{L^1(\Omega)}
\ \lesssim\
\sum_{|\alpha|\le m}\,
\big\|\partial^\alpha f
\big\|_{L^1(\Omega)}\,,
$$
with constant depending only on $C_1$ and $m$.
The result follows.
\hfill
\qed
\vskip.2cm

 Need the preceding result for fractional values of $m$.   
 Define 
 for $0\le m\in \NN$ Banach spaces
 $$
 W^{m,1}(\RR^d)
 :=
 \big\{f\in L^1(\RR^d)\, :\,
 \forall |\alpha|\le m,
 \quad
 \partial^\alpha f\ \in \ L^1(
 \RR^d)
 \big\},
 \quad
 \|
 f\|_{ W^{m,1}(\RR^d)}
  :=
 \sum_{|\alpha|\le m} \|\partial^\alpha f \|_{L^1(\RR^d)}.
 $$
  
 The norms of dilated functions satisfy
  \begin{equation*}
\big\|
f\big(r(\cdot)\big)
\big\|_{L^{1}(\RR^d)}
\, =\,
r^{-(d-1) } \
\big\|
f\big\|_{L^{1}(\RR^d)},
\quad
{\rm and},
\quad
\big\|
\partial^\alpha\big[f\big(r(\cdot)\big)\big]
\big\|_{L^{1}(\RR^d)}
\, =\,
r^{|\alpha|}\
r^{-(d-1) } \,
\big\|
f\big\|_{L^{1}(\RR^d)}.
\end{equation*}
They imply   that for  $f\in W^{m,1}(\RR^d)$
and $r>0$ the function $x\mapsto f(rx)$ also belongs to 
$W^{m,1}(\RR^d)$ and there is a constant independent of $r,f$ so that
\begin{equation}
\label{eq:dilation}
\big\|
f\big(r(\cdot)\big)
\big\|_{W^{m,1}(\RR^d)}
\ \le\
C(m)\
r^{m} \ r^{-(d-1) } \
\big\|
f\big\|_{W^{m,1}(\RR^d)}\,.
\end{equation}

 For $m\ge 1/2$ with $m-1/2\in \NN$,
 $
 \big[
 W^{m-1/2,1}(\RR^d)\,,\,
 W^{m+1/2,1}(\RR^d)
 \big]_{1/2}
$
denotes
 the  complex interpolation space

\begin{lemma}
\label{lem:nonsta2}
With $\Omega, f,\phi$ as above and $m\in \NN +1/2$,
$$
\Big|
\int \ e^{i\phi /\eps}\ 
f(x)\ dx
\Big|
\ \le \
C\, \eps^m\,
\big\|
f
\big\|_{
[
 W^{m-1/2,1}(\RR^d)\,,\,
 W^{m+1/2,1}(\RR^d)
]_{1/2}
}
\ .
$$
\end{lemma}

{\bf Proof.}
Follows by interpolation from the cases $m-1/2$
and $m+1/2$ proved in Lemma \ref{lem:nonsta1}.
\qed

 \begin{definition}  Spaces $Y^m$ are defined for $0\le m\in \NN/2$ as follows.
\begin{align*}
Y^m(\RR^d)
\ &:=\
W^{m,1}(\RR^d)
\qquad\qquad\qquad\qquad\qquad\qquad\ \ 
{\rm for}
\quad
0\le m\in \NN\,, 
\cr
Y^m(\RR^d)\ &:=\
\big[W^{m-1/2,1}(\RR^d)\,,\,
 W^{m+1/2,1}(\RR^d)
\big]_{1/2}
\qquad\,
 {\rm for}\quad
 0<m\in \NN +1/2\,.
\end{align*}
\end{definition}

Lemmas \ref{lem:nonsta1} and \ref{lem:nonsta2} assert that for nonnegative $0\le m\in \NN/2$,
\begin{equation}
\label{eq:dan}
\Big|
\int \ e^{i\phi /\eps}\ 
f(x)\ dx
\Big|
\ \le \
C(m)\
 \eps^m\,
\big\|
f
\big\|_{
Y^m(\RR^d)}
\ .
\end{equation}

 \subsubsection{Limit case inequality of stationary phase}

{\bf Definition.}  {\sl A point $\ux$ in an open subset $\Omega\subset\RR^d$
is a {\bf  stationary point} of
 $\phi\in C^\infty(\Omega\,;\,\RR)$ 
when
$\nabla_x\phi(\ux)=0$.  It is  a {\bf nondegenerate
stationaray point}
when  the matrix of second derivatives at $\ux$
 is nonsingular.
 }

 \vskip.2cm
 
 When $\ux$ is a nondegenerate stationary point,
  the map
 $x\mapsto \nabla_x\phi(x)$ has nonsingular jacobian at 
 $\ux$.  It follows that the map is a local diffeomorphism and
 in particular the stationary point is isolated.
 
 Taylor's Theorem shows that 
 $$
 \nabla_x\phi(x)\ =\ 
 {1\over 2}\,
 \nabla^2_x\phi(\ux)\,(x-\ux) \ +\ 
 O(|x-\ux|^2)\,.
 $$
 If  $\omega \subset\subset\Omega$  contains $\ux$ and
 no other stationary point,
   there is a constant $C>0$ 
 so,
 \begin{equation}
 \label{eq:fred}
\forall \ux\in \omega\,,
\qquad
C^{-1}\,|x-\ux|\ \le\
 \big|
 \nabla_x\phi(x)\big|\ \le\ 
 C\,|x-\ux|
 \,.
 \end{equation}

 The following stationary phase indequality in the 
 borderline regularity case follows a proof I learned
 from G. M\'etiver for the case of more regular $f$
 (see Theroem 3.II.1 in \cite{rauch2012}).  
 
\begin{theorem}
\label{thm:ISP}
 Suppose that 
   $\Omega\subset \RR^n$ is bounded and open and $m$ is the smallest integer
   less than or equal to $n/2$.
   For any $c_1>0$  there is a constant 
$C$ so that for all $f\in C^\infty_0(\Omega)$
 and phase functions $\phi\in C^\infty(\Omega\,;\,\RR)$  there 
a point $\ux\in \Omega$ so that for all $x\in \Omega$
$$
\forall x\in \Omega, \quad  \frac1{c_1}\,|x-\ux|\
\  \le\
\big|
\nabla \phi\big|
\ \le \
c_1\,
 |x-\ux|\,,
 \qquad
 {\rm and},
 \qquad
 \forall |\alpha|\le m+1\,,\quad
 \|
 \partial^\alpha \phi \|_{L^\infty(\Omega)}\ \le \ c_1\,.
 $$
 one has for all $0<\eps<1$
 $$
 \bigg|
 \int
 e^{i\phi(x)/\eps}\
 f(x)\ 
 dx
 \bigg|
  \
  \le
   \
  C\,
  \eps^{n/2}\
  \ln(1+\eps^{-1})\
    \Big(
  \big\|
  f\big\|_{L^\infty(\Gamma)  }
  \ +\  
 \big\|
 f
 \big\|_{Y^{n/2}(\Omega)}
 \Big)
 \,.
     $$
\end{theorem}

  \begin{lemma}
  \label{lem:partition}
   There is  a nonnegative $\chi\in C^\infty_0(\RR^d\setminus 0)$
  so that for all $x\ne 0$,
  $
  \sum_{k=-\infty}^\infty
  \chi(2^k\,x)=
  1
  $.
  \end{lemma}

  {\bf Proof of Lemma.}  Choose nonnegative $g\in C^\infty_0(\RR^d\setminus 0)$ so that
  $g\ge 1$  on $\{1\le |x|\le 2\}$.  Define the locally finite sum
  $$
  G(x)
  \ :=\ 
  \sum_{k=-\infty}^\infty g(2^k x)
  \,,
  \qquad
  {\rm so}
  \qquad
   G(2^kx)=G(x)
  \,.
  $$
  Then $G\in C^\infty(\RR^d\setminus 0)$,
 and $G\ge 1$.
  The function $\chi:= g/G$ has the desired properties.\qed
  
  \vskip.2cm

  {\bf Proof of Theorem.}    
  Translating coordinates it suffices to consider
  $\ux=0$.  Choose $\chi$ as in the Lemma \ref{lem:partition}.  Write
  $$
   \int e^{i\phi/\eps} 
 f(x)\ dx
 \ = \
  \sum_{k=-\infty}^\infty \, \int \ \chi(2^k\,x) \ e^{i\phi/\eps} \
 f(x)\ dx
 \ :=\ \sum_{k=-\infty}^\infty\ I(k)
 \,.
 $$
 
 The half sum $\sum_{k<0} \chi(2^kx)$ is a  smooth function on $\RR^d$ that
  vanishes
 on a neighbhorhood of the origin and is identically equal to 1 outside a large ball.
 Inequality \ref{eq:dan} yields
 $$
 \Big|
 \int  \ e^{i\phi/\eps} 
 \Big(
 \sum_{k<0} \chi(2^k\,x)
  \Big)\
 f(x)
\ dx
 \Big|
 \ \le\
 C\
 \eps^{d/2}
 \
 \big\|
 f\big\|_{Y^{d/2}(\RR^d)}
  \,.
 $$
 
 The sum $\sum_{2^k\eps^{1/2}\ge 1}\chi(2^k\,x)$ is a bounded function
 supported in a ball $|x|\le C\eps^{1/2}$ so
  $$
 \Big|
 \int  \ e^{i\phi/\eps} 
 \Big(
 \sum_{2^k\eps^{1/2}\ge 1} \chi(2^k\,x)
  \Big)\
 f(x)
\ dx
 \Big|
 \ \le\
 C\,
 \eps^{d/2}
 \,
 \| f(x)\|_{L^\infty(\Omega)}
 \,.
 $$

 There remains the  sum over $1\le 2^k< \eps^{-1/2}$.
 The change of variable $y=2^k\,x$ yields
$$
I(k) 
\ =\ 
2^{-kd}\ 
\int
\ 
\chi(y)\
e^{i\phi_k(y)/(2^{2k}\eps)}
\
f(2^{-k}y)
\ 
dy\,,
\qquad
\phi_k(y)
\ :=\ 
2^{2k}\,
\phi(2^{-k}\,y)
\,.
$$
 It follows from \eqref{eq:fred}  that there is a constant $C_1>0$ so that
 on the support of 
 $\chi$,  
 $$
 C_1^{-1}\le \big|\nabla\phi_k\big|\le C_1\,.
 $$
 In addition 
   there is are constants
$C(\alpha)$  independent of $k\ge 0$ so that 
$
|\partial^\alpha\phi_k |
\ \le\
C(\alpha)
$.

Inequality \eqref{eq:dan} shows
 that there is a constant
independent of $k\ge 0$ so that
$$
\Big|
\int
\ 
\chi(y)\
e^{i\phi_k(y)/(2^{2k}\eps)}
\
f(2^{-k}y)
\ 
dy
\Big|
\ \le \ 
C\,
\big(
2^{2k}\,\eps\big)^{d/2}\
\big\|
f
\big\|_{Y^{d/2}(\Omega)}
\ =\
C\
2^{kd}\
\big\|
f
\big\|_{Y^{d/2}(\RR^d)}
 \,.
 $$
Used to estimate $I(k)$, 
the powers of $2^{\pm kd}$ cancel yielding
$$
\sum_{1\le 2^k<\eps^{-1/2}}
|I(k)|
\ \le \ 
C\,\eps^{d/2}\,
\sum_{1\le 2^k<\eps^{-1/2}}
\
\big\|
f\big\|_{Y^{d/2}(\RR^d)}
\,.
$$
The number of summands is $\lesssim \ln(1+\eps^{-1})$.
This completes the proof.
\hfill
\qed
\vskip.2cm

\bibliographystyle{plain}

\end{document}